\documentclass[12pt]{amsart}%
\usepackage{amsmath}
\usepackage{graphicx}
\usepackage{amscd}
\usepackage{geometry}
\usepackage{amsfonts}
\usepackage{amssymb}%
\newtheorem{theorem}{Theorem}[section]
\theoremstyle{plain}

\newtheorem{lemma}[theorem]{Lemma}

\newtheorem{remark}{Remark}

\numberwithin{equation}{section}

\begin{document}
\title[Spherical alterations of handles ]{Spherical alterations of handles: embedding the manifold Plus Construction}
\author{C. R. Guilbault}
\address{Department of Mathematical Sciences, University of Wisconsin-Milwaukee,
Milwaukee, Wisconsin 53201}
\email{craigg@uwm.edu}
\author{F. C. Tinsley}
\address{Department of Mathematics, The Colorado College, Colorado Springs, Colorado 80903}
\email{ftinsley@coloradocollege.edu}
\thanks{This project was aided by a Simons Foundation Collaboration Grant awarded to
the first author. }
\date{January10, 2012}
\subjclass{Primary 57N15, 57Q12; Secondary 57R65, 57Q10}
\keywords{spherical alteration, perfect group, plus construction, generalized plus construction}

\begin{abstract}
Quillen's famous plus construction plays an important role in many aspects of
manifold topology. In our own work on ends of open manifolds, an ability to
embed cobordisms provided by the plus construction into the manifolds being
studied was a key to completing the main structure theorem \cite{GT2}. In this
paper we develop a `spherical modification' trick which allows for a
constructive approach to obtaining those embeddings. More importantly, this
approach can be used to obtain more general embedding results. In this paper
we develop generalizations of the plus construction (together with the
corresponding group-theoretic notions) and show how those cobordisms can be
embedded in manifolds satisfying appropriate fundamental group properties.
Results obtained here are motivated by, and play an important role in, our
ongoing study of noncompact manifolds \cite{GT3}.

\end{abstract}
\maketitle

\section{Introduction}

In this paper we develop a procedure, called \textquotedblleft spherical
alteration\textquotedblright, for modifying handle decompositions of manifolds
in ways that permit useful applications. The strategy is geometrically quite
simple, but at the same time more drastic than the traditional techniques of
handle slides, introductions and cancellations of complementary handle pairs,
and the carving out and inserting of existing handles. In order to obtain the
intended applications, each alteration of a handle is accompanied by
associated alterations of related submanifolds. Taken together, these moves
constitute the process of spherical alteration. Since there are several
variables involved, a full description of the procedure is a bit
technical---we save that for Section \ref{Section: Spherical alteration}. In
some sense, our main result is more a technique than a specific theorem;
nevertheless, several concrete applications of that technique are provided.
The prototypical application is a constructive proof of the following theorem,
which was a key ingredient in the main result \cite{GT2}.

\begin{theorem}
[Embedded Manifold Plus Construction]\label{IntroTh: EPC}Let $R$ be a
connected manifold of dimension $\geq6$, $B$ be a closed component of
$\partial R$, and
\[
K\subseteq\ker\left(  \pi_{1}\left(  B\right)  \rightarrow\pi_{1}\left(
R\right)  \right)
\]
a perfect group that is the normal closure in $\pi_{1}\left(  B\right)  $ of a
finite set of elements. Then there exists an embedding of a plus cobordism
$\left(  W,A,B\right)  $ into $R$ which is the identity on $B$ and for which
$\ker\left(  \pi_{1}\left(  B\right)  \rightarrow\pi_{1}\left(  W\right)
\right)  =K$.
\end{theorem}

\begin{remark}
\emph{(a) Recall that compact cobordism }$\left(  W,A,B\right)  $\emph{ is a
plus cobordism if }$A\hookrightarrow W$\emph{ is a }simple\emph{ homotopy
equivalence. A detailed discussion of plus cobordisms and the Manifold Plus
Construction can be found in Section \ref{Sec: EMPC}.}

\noindent\emph{(b) As an added bonus, our proof of Theorem \ref{IntroTh: EPC}
provides a new twist on the existence proof for plus cobordisms---an argument
that requires very little discussion of bundles and framings.}
\end{remark}

We will further exhibit the usefulness of the spherical alteration technique
by proving a generalization of Theorem \ref{IntroTh: EPC}. That generalization
is motivated by ongoing work on ends of noncompact manifolds. It and a similar
application of spherical alteration, also presented here, play key roles in
\cite{GT3}. To the best of our knowledge, these latter two applications are
not obtainable by the nonconstructive approach to Theorem \ref{IntroTh: EPC}
used in \cite{GT2}.

\section{Preliminaries\label{Section: Preliminaries}}

In this section we provide brief reviews of several topics and introduce a
good deal of notation to be used later in this paper. Those topics are:

\begin{itemize}
\item intersection numbers between submanifolds,

\item surgering surfaces to disks and 2-spheres,

\item perfect groups and `near perfect' subgroups,

\item basics of handle theory, and

\item unbased spheres as elements of homotopy groups.
\end{itemize}

\noindent Throughout this paper we work in the category of piecewise-linear
manifolds; analogous results in the smooth and topological categories may be
obtained in the usual ways.

\subsection{Intersection numbers}

One of the simplest types of intersection number is defined when $P^{p}$ and
$Q^{q}$ are closed, connected, oriented submanifolds of the interior of an
oriented $\left(  p+q\right)  $-manifold $N$. First arrange that $P^{p}$ and
$Q^{q}$ intersect transversely at a finite set of points $p_{1},p_{2}%
,\cdots,p_{k}$. At each $p_{i}$, the local orientation of $P^{p}$ together
with the local orientation of $Q^{q}$ (in that order) determine a local
orientation for $N$. If that orientation agrees with the global orientation of
$N$, we write $sgn\left(  p_{i}\right)  =1$; otherwise $sgn\left(
p_{i}\right)  =-1$. The $\mathbb{Z}$\emph{-intersection number} is defined by
$\varepsilon_{\mathbb{Z}}\left(  P^{p},Q^{q}\right)  =\sum_{i=1}^{k}sgn\left(
p_{i}\right)  $. This definition depends upon order; by linear algebra
$\varepsilon_{\mathbb{Z}}\left(  Q^{q},P^{p}\right)  =\left(  -1\right)
^{pq}\varepsilon_{\mathbb{Z}}\left(  P^{p},Q^{q}\right)  $.

A more delicate intersection \textquotedblleft number\textquotedblright\ lies
in $\mathbb{Z}\left[  \pi_{1}\left(  N,\ast\right)  \right]  $. Instead of
assuming $N$ is oriented (or even orientable), choose a \emph{local}
orientation of $N$ at $\ast$. Assume that $P^{p}$ and $Q^{q}$ are both
oriented and simply connected, and fix \emph{base paths }$\lambda_{P}$ and
$\mu_{Q}$ in $N$ from $\ast$ to base points $\ast_{P}\in P^{p}$ and $\ast
_{Q}\in Q^{q}$. For each $p_{i}$, choose paths $\rho_{i}$ in $P^{p}$ and
$\sigma_{i}$ in $Q^{q}$ from the respective base points to $p_{i}$. Let
$sgn\left(  p_{i}\right)  =\pm1$, depending on whether the local orientation
at $\ast$ translated along the path $\lambda_{P}\cdot\rho_{i}$ agrees with
orientation at $p_{i}$ induced by the orientations of $P^{p}$ then $Q^{q}$;
then let $g_{i}\in\pi_{1}\left(  N,\ast\right)  $ correspond to $\lambda
_{P}\cdot\rho_{i}\cdot\sigma_{i}^{-1}\cdot\mu_{Q}^{-1}$. At $p_{i}$ define
$\varepsilon_{\mathbb{Z\pi}_{1}\left(  N,\ast\right)  }\left(  p_{i}\right)
=sgn\left(  p_{i}\right)  g_{i}$. Finally, the $\mathbb{Z\pi}_{1}%
$\emph{-intersection number} is defined by%
\[
\varepsilon_{\mathbb{Z\pi}_{1}\left(  N,\ast\right)  }\left(  P^{p}%
,Q^{q}\right)  =\sum_{i=1}^{k}\varepsilon_{\mathbb{Z\pi}_{1}\left(
N,\ast\right)  }\left(  p_{i}\right)  \in\mathbb{Z}\left[  \pi_{1}\left(
N,\ast\right)  \right]  .
\]
Note that simple connectivity of $P^{p}$ and $Q^{q}$ ensures that
$\varepsilon_{\mathbb{Z\pi}_{1}\left(  N,\ast\right)  }\left(  P^{p}%
,Q^{q}\right)  $ does not depend on the choice of $\rho_{i}$ and $\sigma_{i}$,
however there is some dependence on $\lambda_{P}$ and $\mu_{Q}$. The ordering
of $P^{p}$ and $Q^{q}$ now plays a larger role than it did for $\mathbb{Z}%
$-intersection numbers---a change in order first alters $sgn\left(
p_{i}\right)  $ by a factor of $\left(  -1\right)  ^{pq}\omega_{1}\left(
g_{i}\right)  $, where $\omega_{1}\left(  g_{i}\right)  =1$ if $g_{i}$ is an
orientation preserving loop and $\omega_{1}\left(  g_{i}\right)  =-1$
otherwise; secondly, the loop $\lambda_{P}\cdot\rho_{i}\cdot\sigma_{i}%
^{-1}\cdot\mu_{Q}^{-1}$ is now traversed in the opposite direction, so $g_{i}$
becomes $g_{i}^{-1}$. For us, the key facts related to order are:

\begin{itemize}
\item $\varepsilon_{\mathbb{Z\pi}_{1}\left(  N,\ast\right)  }\left(
P^{p},Q^{q}\right)  =0$ if and only if $\varepsilon_{\mathbb{Z\pi}_{1}\left(
N,\ast\right)  }\left(  Q^{q},P^{p}\right)  =0$, and

\item if $\varepsilon_{\mathbb{Z\pi}_{1}\left(  N,\ast\right)  }\left(
P^{p},Q^{q}\right)  =1$ then $\varepsilon_{\mathbb{Z\pi}_{1}\left(
N,\ast\right)  }\left(  Q^{q},P^{p}\right)  =\pm1$.
\end{itemize}

Sometimes the simple connectivity conditions on $P^{p}$ and $Q^{q}$ can be
relaxed. An important such case occurs when one of the submanifolds, say
$Q^{q}$, is a $1$-sphere; there we salvage `well-definedness'\ by requiring
that $\sigma_{i}$ be the unique arc of $Q^{q}$ running from $q$ to $p_{i}$ in
the orientation preserving direction. Another useful variation occurs when the
fundamental group of $P^{p}$ or $Q^{q}$ includes trivially into the that of
$N$, in which case that submanifold need not be simply connected. Similarly,
if the images of $\pi_{1}\left(  P^{p}\right)  $ and $\pi_{1}\left(
Q^{q}\right)  $ (translated appropriately along $\lambda_{P}$ and $\mu_{Q}$)
lie in an $L\trianglelefteq\pi_{1}\left(  N,\ast\right)  $ then the above
procedure produces a well-defined intersection number in $\mathbb{Z}\left[
\pi_{1}\left(  N,\ast\right)  /L\right]  $.

We will call collections $\left\{  P_{i}^{p}\right\}  _{i=1}^{r}$ and
$\left\{  Q_{i}^{q}\right\}  _{i=1}^{r}$ of closed submanifolds of a $\left(
p+q\right)  $-manifold $N^{n}$ \emph{geometrically dual} if $P_{i}^{p}$
intersects $Q_{i}^{q}$ transversely in a single point for all $i$ and
$P_{i}^{p}\cap Q_{j}^{q}=\varnothing$ for all $i\neq j$. If $N^{n}$ and all of
these submanifolds are oriented, then the collections are \emph{algebraically
dual} \emph{over} $\mathbb{Z}$ if $\varepsilon_{\mathbb{Z}}\left(  P_{i}%
^{p},Q_{j}^{q}\right)  =\pm\delta_{ij}$ for all $1\leq i,j\leq r$. So (given
the necessary orientability requirements), collections which are geometrically
dual are necessarily algebraically dual over $\mathbb{Z}$, but not conversely.

More generally, given the necessary hypotheses and all required choices to
make $\mathbb{Z[\pi}_{1}\left(  N,\ast\right)  ]$-intersection numbers
well-defined, collections $\left\{  P_{i}^{p}\right\}  _{i=1}^{r}$ and
$\left\{  Q_{i}^{q}\right\}  _{i=1}^{r}$ are \emph{algebraically dual}
\emph{over} $\mathbb{Z[\pi}_{1}\left(  N,\ast\right)  ]$ if $\varepsilon
_{\mathbb{Z[\pi}_{1}\left(  N,\ast\right)  ]}\left(  P_{i}^{p},Q_{j}%
^{q}\right)  =\pm\delta_{ij}$ for all $1\leq i,j\leq r$. In reality, we are
usually satisfied if each $\varepsilon_{\mathbb{Z[\pi}_{1}\left(
N,\ast\right)  ]}\left(  P_{i}^{p},Q_{j}^{q}\right)  =\pm g_{i}$ for some
$g_{i}\in\mathbb{\pi}_{1}\left(  N,\ast\right)  $ and $\varepsilon
_{\mathbb{Z[\pi}_{1}\left(  N,\ast\right)  ]}\left(  P_{i}^{p},Q_{j}%
^{q}\right)  =0$ when $i\neq j$. In those cases, we can always arrange the
more rigid requirement by rechoosing some of the base paths. Under appropriate
conditions, the notion of collections being \emph{algebraically dual}
\emph{over} $\mathbb{Z[\pi}_{1}\left(  N,\ast\right)  /L]$ may be defined in a
similar manner.

\subsection{Surgery on surfaces\label{subsection: surgery on surfaces}}

For a compact oriented surface $\Lambda$ with zero or one boundary components,
a \emph{complete set of meridian-longitude pairs} is a collection of pairs of
oriented simple closed curves $\left\{  \left(  m_{j},l_{j}\right)  \right\}
_{j=1}^{k}$ such that collections $\left\{  m_{j}\right\}  _{j=1}^{k}$and
$\left\{  l_{j}\right\}  _{j=1}^{k}$ are geometrically dual and together
generate $H_{1}\left(  \Lambda;\mathbb{Z}\right)  $.

Given such a collection, let $p_{j}$ denote the point of intersection between
$m_{j}$ and $l_{j}$ and choose a set of arcs $\left\{  \tau_{j}\right\}
_{j=1}^{k}$ in $\Lambda$ intersecting only at a common initial point
$\ast_{\Lambda}$ so that each $\tau_{j}$ intersects the collection of simple
closed curves only at its terminal point $p_{j}$; if $\partial\Lambda
\neq\varnothing$ choose $\ast_{\Lambda}\in\partial\Lambda$. Using
$\Upsilon=\cup_{j=1}^{k}\tau_{j}$ as a `base tree', the curves of $\left\{
\left(  m_{j},l_{j}\right)  \right\}  _{j=1}^{k}$ may be viewed as elements of
$\pi_{1}\left(  \Lambda,\ast_{\Lambda}\right)  $. In the case where
$\partial\Lambda\neq\varnothing$ we may---after relabeling, reordering, and
choosing appropriate orientations on the simple closed curves and on
$\partial\Lambda$---assume that $\partial\Lambda=\prod_{j=1}^{k}m_{j}%
^{-1}l_{j}^{-1}m_{j}l_{j}$ in $\pi_{1}\left(  \Lambda,\ast_{\Lambda}\right)  $.

\begin{remark}
\emph{Since }$\Lambda$\emph{ is not presumed to bound or be embedded in a }%
$3$\emph{-manifold, common distinctions between longitude and meridian (or
neither) are nonexistent here; a given curve could play either role, depending
upon the setup. Nevertheless, the informal use of this terminology will be
convenient for discussing certain curves and collections of curves.}
\end{remark}

Suppose now that $\Lambda$, with zero or one boundary components and a
complete set $\left\{  \left(  m_{j},l_{j}\right)  \right\}  _{j=1}^{k}$ of
meridian-longitude pairs, is embedded in an $n$-manifold $N^{n}$ ($n\geq5$)
and that each $m_{j}$ is homotopically trivial in $N^{n}$. Then we may surger
$\Lambda$ to a $2$-sphere or $2$-disk in the following manner:

\begin{itemize}
\item for convenience, choose a collection $\left\{  \overline{m}_{j}\right\}
_{j=1}^{k}$ of simple closed curves in $\Lambda$ where each $\overline{m}_{j}$
is parallel to $m_{j}$ and disjoint from $\Upsilon.$ Do this so that $\left\{
\overline{m}_{j}\right\}  _{j=1}^{k}$ is geometrically dual to $\left\{
l_{j}\right\}  _{j=1}^{k}$

\item let $\left\{  D_{j}\right\}  _{j=1}^{k}$ be a collection of pairwise
disjoint $2$-disks embedded in $N^{n}$ with $\Lambda\cap D_{j}=\partial
D_{j}=\overline{m}_{j}$,

\item for each $j$, let $A_{j}$ be a small annular neighborhood of
$\overline{m}_{j}$ in $\Lambda$ with boundary curves $\overline{m}%
_{j}^{\hphantom{.}-}$ and $\overline{m}_{j}^{\hphantom{.}+}$,

\item for each $j$, let $D_{j}^{-}$ and $D_{j}^{+}$ be disks parallel to
$D_{j}$ having $\overline{m}_{j}^{\hphantom{.}-}$ and $\overline{m}%
_{j}^{\hphantom{.}+}$ as boundaries,

\item let $\Lambda^{\ast}$ be the $2$-sphere or $2$-disk obtained by removing
the interiors of the $A_{j}$ from $\Lambda$ and sewing in $D_{j}^{-}$ and
$D_{j}^{+}$.
\end{itemize}

If $\Lambda$ has a preferred orientation, there is a corresponding orientation
of $\Lambda^{\ast}$ where the two agree on $\Lambda-\cup A_{j}$. Under that
orientation of $\Lambda^{\ast}$ disk pairs $D_{j}^{\pm}$ inherit opposite
orientations when compared by projecting onto $D_{j}$. Suppose $Q^{n-2}$ is a
closed oriented submanifold of $N^{n}$ intersecting $\Lambda$ transversely in
finitely many points. By applying a small isotopy if necessary, we may assume
none of those intersection points is contained in $\cup A_{j}$. Adjust the
$D_{j}$ (rel boundary) so they also intersect $Q^{n-2}$ transversely.
Corresponding to each $p\in D_{j}\cap Q^{n-2}$ there are points $p^{-}\in
D_{j}^{-}\cap Q^{n-2}$ and $p^{+}\in D_{j}^{+}\cap Q^{n-2}$. Thus
$\Lambda^{\ast}\cap Q^{n-2}$ consists of the points of $\Lambda\cap Q^{n-2}$
together with one pair of points $\left\{  p^{-},p^{+}\right\}  $ for each
point $p$ of a $D_{j}\cap Q^{n-2}$. If $N^{n}$ is oriented, it is clear that
$\varepsilon_{\mathbb{Z}}\left(  p^{-}\right)  =-\varepsilon_{\mathbb{Z}%
}\left(  p^{+}\right)  $ for each of those pairs, so when $\Lambda$ is closed,
we have $\varepsilon_{\mathbb{Z}}\left(  \Lambda^{\ast},Q^{n-2}\right)
=\varepsilon_{\mathbb{Z}}\left(  \Lambda,Q^{n-2}\right)  $.

\begin{figure}
[ptb]
\begin{center}
\includegraphics[
height=2.4693in,
width=3.6292in
]%
{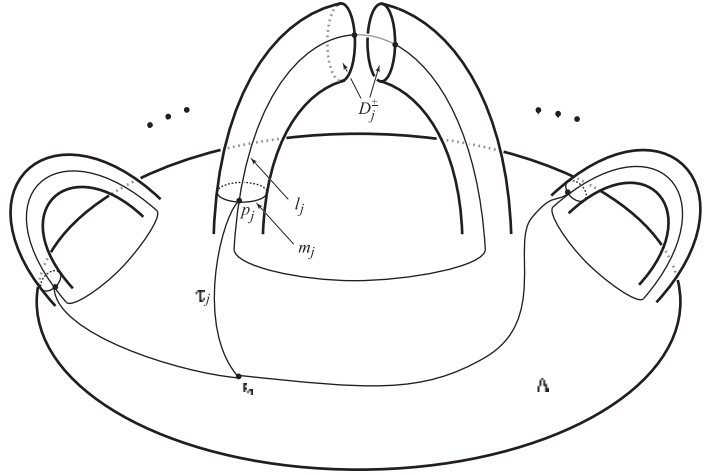}%
\caption{Surgering a surface}%
\label{Fig1}%
\end{center}
\end{figure}

For $\varepsilon_{\mathbb{Z\pi}_{1}\left(  N^{n},\ast\right)  }\left(
\Lambda^{\ast},Q^{n-2}\right)  $ the situation is more complicated. In order
to compare contributions of points $p^{-}$ and $p^{+}$, assume the necessary
setup discussed in the previous subsection: base points $\ast$, $\ast_{Q}$,
$\ast_{\Lambda^{\ast}}=\ast_{\Lambda}$ and corresponding base paths\emph{
}$\mu_{Q}$ and $\lambda_{\Lambda^{\ast}}$, and a local orientation at $\ast$;
assume also that $Q^{n-2}$ is simply connected. Determination of
$\varepsilon_{\mathbb{Z\pi}_{1}\left(  N^{n},\ast\right)  }\left(
p^{-}\right)  $ and $\varepsilon_{\mathbb{Z\pi}_{1}\left(  N^{n},\ast\right)
}\left(  p^{+}\right)  $ require paths $\rho^{-}$ and $\rho^{+}$ in
$\Lambda^{\ast}$ from $\ast_{\Lambda^{\ast}}$ to $p^{-}$ and $p^{+}$. Let
$\rho^{-}$ be the path in $\Lambda^{\ast}$ that follows $\tau_{j}%
\subseteq\Upsilon$ to the longitude $l_{j}$, travels along the unique arc
$l_{j}^{-}\subseteq l_{j}$ that arrives at $D_{j}^{-}$ without leaving
$\Lambda^{\ast}$, and then travels through $D_{j}^{-}$ to $p^{-}$. Choose
$\rho^{+}$ similarly, noting that $l_{j}^{+}$ goes the opposite way around
$l_{j}$. Lastly, choose paths $\sigma^{-}$ and $\sigma^{+}$ in $Q^{n-2}$ from
$\ast_{Q}$ to $p^{-}$ and $p_{+}$; these can be chosen identical except in a
small neighborhood of $p$. If we write $\varepsilon_{\mathbb{Z\pi}_{1}\left(
N,\ast\right)  }\left(  p^{-}\right)  =sgn\left(  p^{-}\right)  g$ and
$\varepsilon_{\mathbb{Z\pi}_{1}\left(  N,\ast\right)  }\left(  p^{+}\right)
=sgn\left(  p^{+}\right)  h$ then $g$ is represented by $\lambda
_{\Lambda^{\ast}}\cdot\rho^{-}\cdot\left(  \sigma^{-}\right)  ^{-1}\cdot
\mu_{Q}^{-1}$ and $h$ by $\lambda_{\Lambda^{\ast}}\cdot\rho^{+}\cdot\left(
\sigma^{+}\right)  ^{-1}\cdot\mu_{Q}^{-1}$. It is easy to check that
$h=l_{j}g$ in $\mathbb{\pi}_{1}\left(  N,\ast\right)  $ and that $sgn\left(
p^{-}\right)  =-\omega_{1}\left(  l_{j}\right)  sgn\left(  p^{+}\right)  $.
See Figure \ref{Fig1}.%
So, together this pair of points contributes $\pm\left(  1-\omega_{1}\left(
l_{j}\right)  l_{j}\right)  g$ to $\varepsilon_{\mathbb{Z\pi}_{1}\left(
N^{n},\ast\right)  }\left(  \Lambda^{\ast},Q^{n-2}\right)  $. For later
reference, we record the following lemma which follows immediately from the
above observations.

\begin{lemma}
\label{Lemma: Surgery on surfaces}Let $\left\{  \left(  m_{j},l_{j}\right)
\right\}  _{j=1}^{k}$ a complete set of meridian-longitude pairs for a closed
oriented surface $\Lambda$ in the interior of an $n$-manifold $N^{n}$ and let
$Q^{n-2}$ be a closed simply connected oriented $\left(  n-2\right)
$-manifold also lying in $\operatorname{int}N^{n}$ and intersecting $\Lambda$
transversely. Assume that each of the meridianal curves $m_{j}$ contracts in
$N^{n}$ and let $\Lambda^{\ast}$ be a $2$-sphere obtained by surgering
$\Lambda$ along a collection of parallel curves; do this in such a way that
$\Lambda^{\ast}$ and $Q^{n-2}$ intersect transversely. Choose base points
$\ast$, $\ast_{\Lambda}=\ast_{\Lambda^{\ast}}$, and $\ast_{Q}$, base
paths\emph{ }$\lambda_{\Lambda}$ and $\mu_{Q}$ from $\ast$ to $\ast_{\Lambda}$
and $\ast_{Q}$, respectively and a local orientation of $N^{n}$ at $\ast$. Then

\begin{enumerate}
\item If each longitudinal curve $l_{i}$ also contracts in $N^{n}$, then
$\varepsilon_{\mathbb{Z\pi}_{1}\left(  N^{n},\ast\right)  }\left(
\Lambda,\allowbreak Q^{n-2}\right)  $ is well-defined and equal to
$\varepsilon_{\mathbb{Z\pi}_{1}\left(  N^{n},\ast\right)  }\left(
\Lambda^{\ast},Q^{n-2}\right)  $.

\item If $L$ is a normal subgroup of $\pi_{1}\left(  N^{n},\ast\right)  $ with
$\omega_{1}\left(  L\right)  \equiv1$ and each $l_{i}$ represents an element
of $L$, then $\varepsilon_{\mathbb{Z[\pi}_{1}\left(  N^{n},\ast\right)
/L]}\left(  \Lambda,Q^{n-2}\right)  $ is well-defined and equal to
$\varepsilon_{\mathbb{Z[\pi}_{1}\left(  N^{n},\ast\right)  /L]}\left(
\Lambda^{\ast},Q^{n-2}\right)  $.
\end{enumerate}
\end{lemma}

\subsection{Perfect groups and nearly perfect
subgroups\label{Subsection: perfect and nearly perfect groups}}

The \emph{commutator subgroup }of a group\emph{ }$G$, denoted $\left[
G,G\right]  $, is the subgroup generated by all commutator elements $\left[
m,l\right]  =m^{-1}l^{-1}ml$ for $l,m\in G$. It is standard knowledge that
$\left[  G,G\right]  $ is a normal and that $G/\left[  G,G\right]  $ is
abelian; in fact, $\left[  G,G\right]  $ is the smallest subgroup of $G$ with
abelian quotient. If $G=\left[  G,G\right]  $, or equivalently $G/\left[
G,G\right]  $ is trivial, we say $G$ is \emph{perfect}. In this paper, we are
interested in topological implications of these concepts. If $G=\pi_{1}\left(
X,x\right)  $ and $\alpha=\prod_{j=1}^{k}\left[  m_{j},l_{j}\right]
\in\left[  G,G\right]  $ then there exists a (mapped in) compact orientable
surface $\Lambda_{\alpha}$ with boundary corresponding to $\alpha$ and a base
tree for which a complete set of meridian-longitude pairs has the form
$\left\{  \left(  m_{j},l_{j}\right)  \right\}  _{j=1}^{k}$. If $\alpha\in K$,
where $K$ is a perfect subgroup of $G$, we may arrange that all of the $m_{j}$
and $l_{j}$ are elements of $K$; this is a key property of perfect subgroups
of $\pi_{1}\left(  X,x\right)  $.

Next we generalize the notion of \textquotedblleft
perfectness\textquotedblright\ for subgroups of $G$. Suppose $K\leq L\leq G$
where $K$ and $L$ are normal in $G$. Then $\left[  K,L\right]  $ is the
subgroup of $G$ generated by all elements of the form $\left[  m,l\right]  $
where $m\in K$ and $l\in L$. It is easy to see that $\left[  K,L\right]
=\left[  L,K\right]  $ and $\left[  K,L\right]  \leq K$. We say that $K$ is
\emph{strongly }$L$\emph{-perfect} if $K\leq\left[  K,L\right]  $. Clearly,
\thinspace$K$ is perfect if and only if it is strongly $K$-perfect; more
generally, the smaller the subgroup $L$ containing $K$, the closer a strongly
$L$-perfect group is to being perfect. When $G=\pi_{1}\left(  X,x\right)  $,
for each element $\alpha$ of a strongly $L$-perfect group $K$, there exists a
(mapped in) compact orientable surface $\Lambda_{\alpha}$ with boundary
corresponding to $\alpha$ and a base tree for which a complete set of
meridian-longitude pairs has the form $\left\{  \left(  m_{j},l_{j}\right)
\right\}  _{j=1}^{k}$ where each $m_{j}\in K$ and $l_{j}\in L$.

\begin{remark}
\emph{We have reserved the term `}$L$\emph{-perfect' (as compared to `strongly
}$L$\emph{-perfect') for the case }$K\leq\left[  L,L\right]  $\emph{, a weaker
condition that is developed in \cite{GT2} but is not used here.}
\end{remark}

\subsection{Basic handle theory\label{subsection: handle theory}}

Let $N^{n}$ be an $n$-manifold, $B$ a component of $\partial N^{n}$ and $J$ a
subset of $B$ homeomorphic to $S^{k-1}\times D^{n-k}$. The act of attaching a
$k$\emph{-handle} (or a \emph{handle of index }$k$) to $N^{n}$ along $J$ is
the creation of an adjunction space $N^{n}\cup_{f}D^{n}$, where $D^{n}$ is
viewed as $D^{k}\times D^{n-k}$ and $f:S^{k-1}\times D^{n-k}\rightarrow J$ is
a homeomorphism. We denote the adjunction space by $N^{n}\cup h^{k}$; here
$h^{k}$ denotes the image of $D^{n}$ under the quotient map $q:N^{n}\sqcup
D^{n}\rightarrow N^{n}\cup_{f}D^{n}$. We call $J$ the \emph{attaching tube} of
$h^{k}$ and $\alpha^{k-1}=q(S^{k-1}\times\left\{  0\right\}  )$ the
\emph{attaching sphere}. We call $e^{k}=q\left(  D^{k}\times\left\{
0\right\}  \right)  $ the \emph{core} and $q\left(  \left\{  0\right\}  \times
D^{n-k}\right)  $ the \emph{cocore} of $h^{k}$; the boundary of the cocore,
$\beta^{n-k-1}=q\left(  \left\{  0\right\}  \times S^{n-k-1}\right)  $, is the
\emph{belt sphere} and $q\left(  D^{k}\times S^{n-k-1}\right)  $ is the
\emph{belt tube} of $h^{k}$. We refer to the boundary component of $N^{n}\cup
h^{k}$ consisting of $B-J$ and the belt tube of $h^{k}$ informally as the
\emph{right-hand boundary component. }The homeomorphism $f:S^{k-1}\times
D^{n-k}\rightarrow J$ is called the \emph{framing }of $h^{k}$; it affects the
homeomorphism type of $N^{n}\cup h^{k}$. More explicitly, if $f_{1}%
,f_{2}:S^{k-1}\times D^{n-k}\rightarrow J$ is a pair of framings, then the
resulting manifolds $N^{n}\cup h_{1}^{k}$ and $N^{n}\cup h_{2}^{k}$ need not
be homeomorphic.

The pair $(N^{n}\cup h^{k},N^{n})$ is homotopy equivalent to $(N^{n}\cup
e^{k},N^{n})$ in the obvious manner. Thus, if $W^{n}=N^{n}\cup h_{1}\cup
\cdots\cup h_{r}$ is obtained from $N^{n}$ by successive attachment of handles
of non-decreasing index, each to the right-hand boundary of the preceding
space, then $\left(  W^{n},N^{n}\right)  $ is homotopy equivalent to a
relative CW complex $\left(  K,N^{n}\right)  $ with one $j$-cell for each
$j$-handle. A useful relationship between these spaces is the equivalence of
the $\mathbb{Z}$-incidence number $\varepsilon_{\mathbb{Z}}\left(
e^{j+1},e^{j}\right)  $ for a pair of cells $e^{j+1}$ and $e^{j}$ of $K$ and
the $\mathbb{Z}$-intersection number $\varepsilon_{\mathbb{Z}}\left(
\alpha^{j},\beta^{n-j-1}\right)  $ of the attaching sphere and the belt sphere
of corresponding $\left(  j+1\right)  $- and $j$-handles, and similarly for
$\mathbb{Z\pi}_{1}$-incidence and intersection numbers. (Both of these
observations require a careful setup involving base points, base paths,
orientations, etc. and some special care in dealing with those cases where the
attaching or belt sphere is not simply connected. The reader is referred to
\cite{RS} for details.) The upshot of all this is that intersection numbers,
employed appropriately, allow one to calculate algebraic data such as
$\mathbb{Z}$-homology, $\mathbb{Z\pi}_{1}$-homology, and Whitehead torsion for
manifolds constructed through the addition of handles.

\subsection{Unbased $k$-spheres as elements of $\pi_{k}\left(  N,\ast\right)
\label{Subsection: unbased k-spheres}$}

Let $\Sigma_{1}$ and $\Sigma_{2}$ be embedded oriented $k$-spheres ($k\geq2$)
in a connected manifold $N$, $\ast\in N$, $\ast_{1}\in\Sigma_{1}$, and
$\ast_{2}\in\Sigma_{2}$. Let $S^{k}$ be the standard $k$-sphere with
$\mathbf{e}_{1}=\left(  1,0,\cdots,0\right)  $ the canonical base point. In
order to view $\Sigma_{1}$ as an element of $\pi_{k}\left(  N,\ast\right)  $,
choose a path $\lambda$ from $\ast$ to $\ast_{1}$. Now define a map of
$(S^{k},\mathbf{e}_{1})$ into $\left(  N,\ast\right)  $ as follows. Let
$D^{k}\subseteq S^{k}$ be a small $k$-disk centered at $\mathbf{e}_{1}$. Send
$\mathbf{e}_{1}$ to $\ast$ and the radial lines of $D^{k}$ emanating from
$\mathbf{e}_{1}$ each onto $\lambda$; then send $S^{k}-D^{k}$ homeomorphically
onto $\Sigma_{1}-\left\{  \ast_{1}\right\}  $ in an orientation preserving
manner. We denote the corresponding element of a $\pi_{k}\left(
N,\ast\right)  $ by $\left[  \lambda\Sigma_{1}\right]  $. Stated differently,
$\left[  \lambda\Sigma_{1}\right]  $ is the image of $\left[  \Sigma
_{1}\right]  \in\pi_{k}\left(  N,\ast_{1}\right)  $ under the change of base
points isomorphism induced by $\lambda$. (See, for example, \cite[\S 4.1]{Ha}.)

\begin{remark}
\emph{If }$\ast_{1}=\ast$\emph{, then }$\lambda$\emph{ is a loop and we
represent }$\left[  \lambda\Sigma_{1}\right]  $\emph{ by }$\lambda\left[
\Sigma_{1}\right]  $\emph{, the image of }$\left[  \Sigma_{1}\right]  $\emph{
when acted upon by }$\lambda$\emph{ under well-known action of }$\pi
_{1}\left(  N,\ast\right)  $\emph{ on }$\pi_{k}\left(  N,\ast\right)  $\emph{.
Since }$\pi_{k}\left(  N,\ast\right)  $\emph{ is abelian, this action may be
extended in the obvious way to an action of }$Z\left[  \pi_{1}\left(
N,\ast\right)  \right]  $\emph{ on }$\pi_{k}\left(  N,\ast\right)  $\emph{.
Again see \cite[\S 4.1]{Ha}.}
\end{remark}

Returning to the original setup, if $\lambda^{\prime}$ is another path from
$\ast$ to $\ast_{1}$, then $\left[  \lambda\Sigma_{1}\right]  $ and $\left[
\lambda^{\prime}\Sigma_{1}\right]  $ need not be equal; it is easy to see that
$\left[  \lambda^{\prime}\Sigma_{1}\right]  =\left(  \lambda^{\prime}%
\cdot\lambda^{-1}\right)  \left[  \lambda\Sigma_{1}\right]  $.

Now suppose $\xi$ is a path from a point $\ast^{\prime}\in\Sigma_{1}$ to
$\ast_{2}$. By a strategy similar to the above, we may obtain a map of $S^{k}$
into $N$ sending a slightly shrunken lower hemisphere onto $\Sigma_{1}$, a
slightly shrunken upper hemisphere onto $\Sigma_{2}$ (both in orientation
preserving manners), and taking a product neighborhood of $S^{k-1}$ onto $\xi
$. When the codimension is sufficiently high we may obtain an embedded version
of the above, denoted $\Sigma_{1}\#_{\xi}\Sigma_{2}$, it consists of punctured
copies of $\Sigma_{1}$ and $\Sigma_{2}$ joined by a `tube'\ following $\xi$.
In either case, we express the corresponding element of $\pi_{k}\left(
N,\ast_{1}\right)  $ by $\left[  \Sigma_{1}\#_{\xi}\Sigma_{2}\right]  $.
Returning to our preferred base point we have $\left[  \lambda\left(
\Sigma_{1}\#_{\xi}\Sigma_{2}\right)  \right]  \in\pi_{k}\left(  N,\ast\right)
$, an element that may be expressed as an algebraic sum as follows: choose a
path $\sigma$ in $\Sigma_{1}$ from $\ast_{1}$ to $\ast^{\prime}$, then
$\left[  \lambda\left(  \Sigma_{1}\#_{\xi}\Sigma_{2}\right)  \right]  =\left[
\lambda\Sigma_{1}\right]  +\left[  \left(  \lambda\cdot\sigma\cdot\xi\right)
\Sigma_{2}\right]  .$

\section{Spherical alteration of a handle\label{Section: Spherical alteration}%
}

In this section we give a precise formulation of spherical alteration and
prove the corresponding technical lemma. Although the technique can be applied
to handles of various indices (usually with some restrictions on codimension),
all of our current applications involve alterations of $2$-handles. For that
reason, we restrict attention to $2$-handles and invite the reader to consider
possible applications of higher index alterations.

\begin{figure}
[ptb]
\begin{center}
\includegraphics[
height=3.0851in,
width=4.1289in
]%
{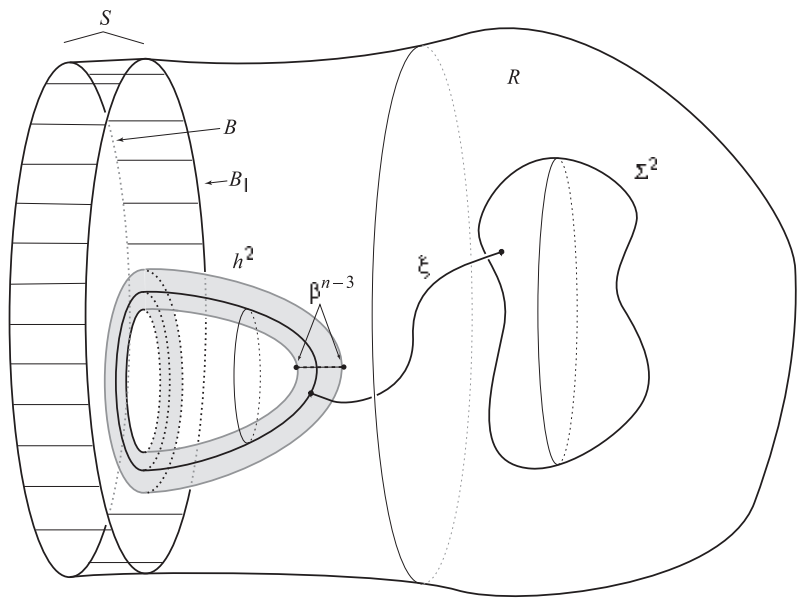}%
\caption{Setup for a spherical alteration}%
\label{Fig2}%
\end{center}
\end{figure}

Let $R$ be an $n$-manifold of dimension at least $6$, $B$ a codimension $0$
submanifold of $\partial R$, $S\approx B\times\left[  0,1\right]  $ a collar
neighborhood of $B$ in $R$, and $h^{2}$ a $2$-handle attached to the interior
boundary component $B_{1}$ of $S$ and lying in $\overline{R-S}$. Let $T=S\cup
h^{2}$, $B_{2}=\partial T-B$, and $e^{2}$ the core of $h^{2}$. In addition,
let $\Sigma^{2}$ be an oriented $2$-sphere embedded in the interior of $R-T$
and $\xi$ an arc in $R$ from a point $p\in e^{2}$ to $q\in\Sigma^{2}$,
intersecting $e^{2}\cup\Sigma^{2}$ at no other points. See Figure \ref{Fig2}%
.

The \emph{spherical alteration of }$h^{2}$ \emph{over }$\Sigma^{2}$
\emph{along }$\xi$ is another $2$-handle in $R$ with the same attaching tube
as $h^{2}$, but with a core $e^{2}\#_{\xi}\Sigma^{2}$ (the connected sum of
$e^{2}$ with $\Sigma^{2}$ along a tube\ contained in a regular neighborhood of
$\xi$). We will denote this new $2$-handle by $h^{2}(\xi,\Sigma^{2})$. An
orientation on $e^{2}$ (induced by a preferred characteristic map for $h^{2}$)
and on $\Sigma^{2}$ are necessary to define $e^{2}\#_{\xi}\Sigma^{2}$; the
connecting tube must be chosen to respect these orientations. More precisely,
let $E\subseteq E^{\prime}\subseteq e^{2}$ and $F\subseteq\Sigma^{2}$ be small
$2$-disks centered at $p$ and $q$, and $Z_{\xi}$ an embedded copy of
$S^{1}\times\left[  0,1\right]  $ contained in a regular neighborhood of $\xi$
such that $Z_{\xi}\cap e^{2}=\partial E$ and $Z_{\xi}\cap\Sigma^{2}=\partial
F$ (at opposite ends of $Z_{\xi}$). Then $e^{2}\#_{\xi}\Sigma^{2}=\left(
e^{2}-\overset{\circ}{E}\right)  \cup Z_{\xi}\cup\left(  \Sigma^{2}%
-\overset{\circ}{F}\right)  $. If the orientation on $e^{2}-\overset{\circ}%
{E}$ translated along $Z_{\xi}$ does not match the orientation on $\Sigma
^{2}-\overset{\circ}{F}$, use the extra codimension to rechoose $Z_{\xi}$ with
a twist so that the orientations match.

\begin{figure}
[ptb]
\begin{center}
\includegraphics[
height=2.8144in,
width=5.0275in
]%
{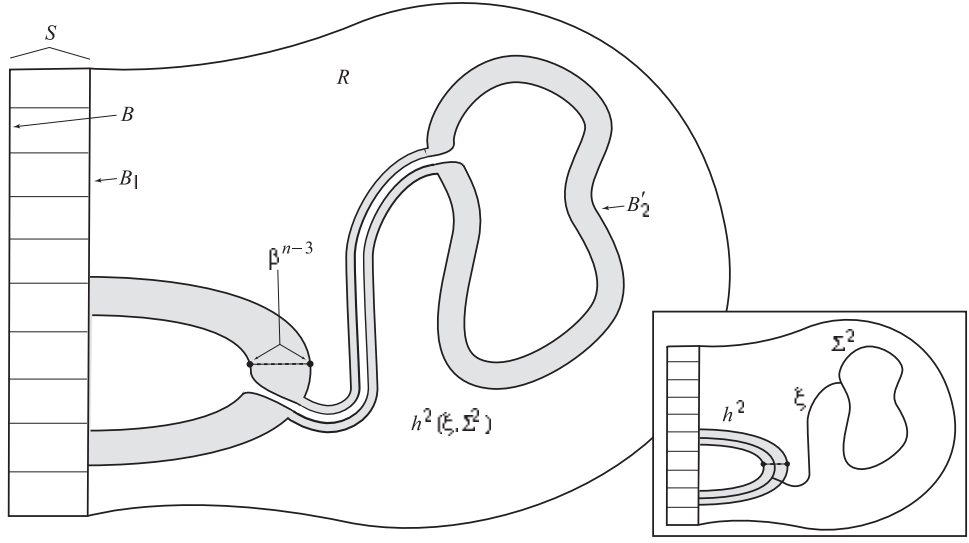}%
\caption{Schematic of a spherical alteration}%
\label{Fig3}%
\end{center}
\end{figure}

Use a parameterization of $h^{2}$ as $D^{2}\times D^{n-2}$ to identify
$P,\mathring{P}\subseteq h^{2}$ corresponding to $E^{\prime}\times D^{n-2}$
and $\overset{\circ}{E^{\prime}}\times D^{n-2}$, respectively; we will refer
to $P$ as the \emph{alteration region} in $h^{2}$. Let $\widehat{h}%
=h^{2}-\mathring{P}$ and choose a relative regular neighborhood $N$ of the
$2$-disk $\left(  E^{\prime}-\overset{\circ}{E}\right)  \cup Z_{\xi}%
\cup\left(  \Sigma^{2}-\overset{\circ}{F}\right)  $ in $\overline{R-\left(
S\cup\widehat{h}\right)  \text{ }}$so that $N$ intersects the boundary
precisely in $P-\mathring{P}\approx S^{1}\times D^{n-2}$. Define $h^{2}%
(\xi,\Sigma^{2})$ to be $\widehat{h}\cup N$, which is a regular neighborhood
of $e^{2}\#_{\xi}\Sigma^{2}$ in $\overline{R-S}$, and thus a $2$-handle.
Clearly $h^{2}(\xi,\Sigma^{2})$ has the same attaching tube as $h^{2}$
(possibly with different framing). We also identify a common belt sphere for
$h^{2}$ and $h^{2}(\xi,\Sigma^{2})$ lying just outside the alteration region:
let $z\in D^{2}$ correspond to a point of $e^{2}$ lying just outside
$E^{\prime}$ and let $\beta^{n-3}$ be the $\left(  n-3\right)  $-sphere
corresponding to $\partial\left(  z\times D^{n-2}\right)  $. Finally, let
$B_{2}^{\prime}$ denote the right-hand boundary of $T^{\prime}=S\cup h^{2}%
(\xi,\Sigma^{2})$. See Figure \ref{Fig3}%
.

Now assume that, in addition to the above, there is a $2$-sphere $\Gamma^{2}$
lying in $B_{2}$ and transverse to $\beta^{n-3}$. If $\Gamma^{2}$ and
$\beta^{n-3}$ intersect in an essential way, then $\Gamma^{2}$ will not lie in
$B_{2}^{\prime}$. Instead, each \textquotedblleft sheet\textquotedblright\ of
$\Gamma^{2}$ that cuts through $\beta^{n-3}$ leaves $B_{2}^{\prime}$ at the
alteration region. We wish to define an alteration of $\Gamma^{2}$ to a
$2$-sphere that lies in $B_{2}^{\prime}$ and intersects $\beta^{n-3}$ in the
same way that $\Gamma^{2}$ does. Let $\Gamma^{2}\cap\beta^{n-3}=\left\{
p_{1},\cdots,p_{k}\right\}  $. Using the product structure of the belt tube,
we may arrange (via an ambient isotopy) that $\Gamma^{2}\cap h^{2}=\left\{
D_{1},\cdots D_{k}\right\}  $ where the $D_{i}$'s are $2$-disks in the belt
tube parallel to $e^{2}$. Remove from each $D_{i}$ the interior of the subdisk
$D_{i}^{\prime}=D_{i}\cap P$, and replace it with an oriented disk
$D_{i}^{\prime\prime}$ which has the same boundary (and induces the same
orientation on that boundary), but lies in the boundary of the regular
neighborhood $N$. Note that $D_{i}^{\prime\prime}$ will be parallel in $N$ to
the core of $N$, with the two having matching orientations if $sgn\left(
p_{i}\right)  =1$ and opposite orientations when $sgn\left(  p_{i}\right)
=-1$. (Orientations are compared by retracting $N$ onto its core.) By general
position, we may assume the $D_{i}^{\prime\prime}$ are disjoint; thus we have
a new $2$-sphere $\Gamma^{2}(\xi,\Sigma^{2})$, called the \emph{corresponding
spherical alteration of }$\Gamma^{2}$ \emph{over }$\Sigma^{2}$ \emph{along
}$\xi$.

Our main lemma equates $\varepsilon_{\mathbb{Z[\pi}_{1}\left(  B_{2}%
,\ast\right)  ]}\left(  \Gamma^{2},\beta^{n-3}\right)  $ with $\varepsilon
_{\mathbb{Z[\pi}_{1}\left(  B_{2}^{\prime},\ast\right)  ]}\left(  \Gamma
^{2}(\xi,\Sigma^{2}),\beta^{n-3}\right)  $ and compares $\Gamma^{2}$ and
$\Gamma^{2}(\xi,\Sigma^{2})$ when viewed as elements of $\pi_{2}\left(
R,\ast\right)  $. Making our assertions precise requires preliminary work.
First, choose $\ast$ to lie in a portion of $B_{1}$ away from $h^{2}$; then
$\ast$ lies in both $B_{2}$ and $B_{2}^{\prime}$. Although $T$ need not be
homeomorphic to $T^{\prime}$ (attaching tubes of the $2$-handles are the same,
but framings may differ), the fundamental groups are canonically
isomorphic---both are obtained by taking the quotient of $\pi_{1}\left(
B_{1},\ast\right)  $ by the normal closure of the common attaching circle for
the $2$-handles. Since $B_{2}\hookrightarrow T$ and $B_{2}^{\prime
}\hookrightarrow T^{\prime}$ induce $\pi_{1}$-isomorphisms, there is a
canonical isomorphism $\phi:\pi_{1}\left(  B_{2},\ast\right)  \rightarrow
\pi_{1}\left(  B_{2}^{\prime},\ast\right)  $ that associates each loop in
$B_{2}$ missing the alteration region\ (a collection that generates the entire
fundamental group) to the identical loop in $B_{2}^{\prime}$. In a casual
sense, $\phi$ may be viewed as an identity map. Let $\widehat{\phi
}:\mathbb{Z[\pi}_{1}\left(  B_{2},\ast\right)  ]\rightarrow\mathbb{Z[\pi}%
_{1}\left(  B_{2}^{\prime},\ast\right)  ]$ be the corresponding
\textquotedblleft identity-like\textquotedblright\ ring isomorphism.

In order to calculate\ $\varepsilon_{\mathbb{Z[\pi}_{1}\left(  B_{2}%
,\ast\right)  ]}\left(  \Gamma^{2},\beta^{n-3}\right)  $ we select some
initial data. Assume that $\left\{  p_{1},\cdots,p_{k}\right\}  $ is non-empty
and choose $p_{1}$ as base point for both $\Gamma^{2}$ and $\beta^{n-3}$.
Choose a local orientation for $B_{2}$ at $\ast$ and a path $\lambda$ in
$B_{2}$ from $\ast$ to $p_{1}$ to serve as base path for both $\Gamma^{2}$ and
$\beta^{n-3}$. The remaining data needed for $\varepsilon_{\mathbb{Z[\pi}%
_{1}\left(  B_{2},\ast\right)  ]}\left(  \Gamma^{2},\beta^{n-3}\right)  $ is a
collection of paths $\rho_{i}$ in $\Gamma^{2}$ and $\sigma_{i}$ in
$\beta^{n-3}$ and from $p_{1}$ to $p_{i}$ for each $i=2,\cdots,k$. ($\rho_{1}$
and $\sigma_{1}$ can be constant paths.). When calculating $\varepsilon
_{\mathbb{Z[\pi}_{1}\left(  B_{2}^{\prime},\ast\right)  ]}\left(  \Gamma
^{2}(\xi,\Sigma^{2}),\beta^{n-3}\right)  $ notice that $\Gamma^{2}(\xi
,\Sigma^{2})\cap\beta^{n-3}=\left\{  p_{1},\cdots,p_{k}\right\}  $ and that
the local orientation at $\ast$ and each of the points and paths just selected
can be chosen to lie simultaneously in $B_{2}$ and $B_{2}^{\prime}$; one
simply avoids the alteration region. By using the same initial data for both,
it is immediate that $\varepsilon_{\mathbb{Z[\pi}_{1}\left(  B_{2}%
,\ast\right)  ]}(\Gamma^{2},\beta^{n-3})$ is the same as $\varepsilon
_{\mathbb{Z[\pi}_{1}\left(  B_{2}^{\prime},\ast\right)  ]}\left(  \Gamma
^{2}(\xi,\Sigma^{2}),\beta^{n-3}\right)  $ --- or more precisely,
$\widehat{\phi}$ takes the former to the latter.

Two more items are needed in preparation for the statement of our main lemma.
First, for each $i=1,\cdots,k$, choose a path $\xi_{i}$ from $p_{i}$ to the
point $q\in\Sigma^{2}$ which travels the short distance through the belt tube
from $p_{i}$ to the alteration region, then runs parallel to $\xi$ and ends at
$q$. It is then clear that each $\xi_{i}$ is homotopic in $R$ (rel endpoints)
to $\sigma_{i}^{-1}\cdot\xi_{1}$. Finally, the inclusion $j:B_{2}%
\hookrightarrow R$ induces a group homomorphism $j_{\#}:\mathbb{\pi}%
_{1}\left(  B_{2},\ast\right)  \rightarrow\mathbb{\pi}_{1}\left(
R,\ast\right)  $ and corresponding ring homomorphism $\widehat{j}%
_{\#}:\mathbb{Z[\pi}_{1}\left(  B_{2},\ast\right)  ]\rightarrow\mathbb{Z[\pi
}_{1}\left(  R,\ast\right)  ]$ which plays a role in the following.

\begin{lemma}
[Spherical Alteration Lemma]\label{alteration lemma}Given the spherical
alteration of $h^{2}$ over $\Sigma^{2}$ along $\xi$ described above, the
corresponding alteration of $\Gamma^{2}$, and all of the base point, path, and
homomorphism data selected in the previous three paragraphs, the following is true:

\begin{enumerate}
\item $\varepsilon_{\mathbb{Z[\pi}_{1}\left(  B_{2},\ast\right)  ]}(\Gamma
^{2},\beta^{n-3})$ is taken to $\varepsilon_{\mathbb{Z[\pi}_{1}\left(
B_{2}^{\prime},\ast\right)  ]}\left(  \Gamma^{2}(\xi,\Sigma^{2}),\beta
^{n-3}\right)  $ by $\widehat{\phi}:\mathbb{Z[\pi}_{1}\left(  B_{2}%
,\ast\right)  ]\allowbreak\rightarrow\allowbreak\mathbb{Z[\pi}_{1}\left(
B_{2}^{\prime},\ast\right)  ]$, and

\item as elements of $\pi_{2}\left(  R,\ast\right)  $, $[\lambda\Gamma^{2}%
(\xi,\Sigma^{2})]$ is equal to%
\[
\left[  \lambda\Gamma^{2}\right]  +\widehat{j}_{\#}(\varepsilon_{\mathbb{Z[\pi
}_{1}\left(  B_{2},\ast\right)  ]}\left(  \Gamma^{2},\beta^{n-3}\right)
\left[  (\lambda\cdot\xi_{1})\Sigma^{2}\right]  ).
\]

\end{enumerate}

\begin{proof}
Item (1) was covered in the lead-up to this lemma.

As for item (2), it is easy to see that the embedding which takes $S^{2}$ onto
$\Gamma^{2}(\xi,\Sigma^{2})$ is homotopic in $R$ to the map indicated by%
\[
\Gamma^{2}\#_{\xi_{1}}(sgn\left(  p_{1}\right)  \Sigma^{2})\#_{\xi_{2}%
}(sgn\left(  p_{2}\right)  \Sigma^{2})\#_{\xi_{3}}\cdots\#_{\xi_{k}}\left(
sgn\left(  p_{k}\right)  \right)  \Sigma^{2}\text{,}%
\]
where a minus sign indicates a reversed orientation. By repeatedly applying
the observations made in Section \ref{Subsection: unbased k-spheres}, we see
that
\begin{align*}
&  \left[  \lambda(\Gamma^{2}\#_{\xi_{1}}(sgn\left(  p_{1}\right)  \Sigma
^{2})\#_{\xi_{2}}(sgn\left(  p_{2}\right)  \Sigma^{2})\#_{\xi_{3}}%
\cdots\#_{\xi_{k}}\left(  sgn\left(  p_{k}\right)  \right)  \Sigma^{2})\right]
\\
&  =\left[  \lambda\Gamma^{2}\right]  +\left[  \left(  \lambda\cdot\xi
_{1}\right)  sgn\left(  p_{1}\right)  \Sigma^{2}\right]  +\sum_{i=2}%
^{k}sgn\left(  p_{i}\right)  \left[  \left(  \lambda\cdot\rho_{i}\cdot\xi
_{i}\right)  \Sigma^{2}\right] \\
&  =\left[  \lambda\Gamma^{2}\right]  +sgn\left(  p_{1}\right)  \left[
\left(  \lambda\cdot\xi_{1}\right)  \Sigma^{2}\right]  +\sum_{i=2}%
^{k}sgn\left(  p_{i}\right)  [(\lambda\cdot\rho_{i}\cdot\left(  \sigma
_{i}^{-1}\cdot\xi_{1}\right)  )\Sigma^{2}]\\
&  =\left[  \lambda\Gamma^{2}\right]  +sgn\left(  p_{1}\right)  \left[
\left(  \lambda\cdot\xi_{1}\right)  \Sigma^{2}\right]  +\sum_{i=2}%
^{k}sgn\left(  p_{i}\right)  [(\lambda\cdot\rho_{i}\cdot\sigma_{i}^{-1}%
\cdot(\lambda^{-1}\cdot\lambda)\cdot\xi_{1})\Sigma^{2}]\\
&  =\left[  \lambda\Gamma^{2}\right]  +sgn\left(  p_{1}\right)  \left[
\left(  \lambda\cdot\xi_{1}\right)  \Sigma^{2}\right]  +\sum_{i=2}%
^{k}sgn\left(  p_{i}\right)  \left(  \lambda\cdot\rho_{i}\cdot\sigma_{i}%
^{-1}\cdot\lambda^{-1}\right)  \left[  \left(  \lambda\cdot\xi_{1}\right)
\Sigma^{2}\right] \\
&  =\left[  \lambda\Gamma^{2}\right]  +\left[  \left(  \lambda\cdot\xi
_{1}\right)  \Sigma^{2}\right]  +\sum_{i=2}^{k}sgn\left(  p_{i}\right)
g_{i}\left[  \left(  \lambda\cdot\xi_{1}\right)  \Sigma^{2}\right] \\
&  =\left[  \lambda\Gamma^{2}\right]  +\left(  sgn\left(  p_{1}\right)
+\sum_{i=2}^{k}sgn\left(  p_{i}\right)  g_{i}\right)  \left[  \left(
\lambda\cdot\xi_{1}\right)  \Sigma^{2}\right]
\end{align*}
where $g_{i}=\lambda\cdot\rho_{i}\cdot\sigma_{i}^{-1}\cdot\lambda^{-1}$ is
precisely the loop used in defining $\varepsilon_{\mathbb{Z[\pi}_{1}\left(
B_{2},\ast\right)  ]}\left(  p_{i}\right)  $ for $i=2,\cdots,k$. By our choice
of base paths, the loop corresponding to $p_{1}$ is null-homotopic, so
$\varepsilon_{\mathbb{Z[\pi}_{1}\left(  B_{2},\ast\right)  ]}\left(
p_{1}\right)  =sgn\left(  p_{1}\right)  $. Thus $\left(  sgn\left(
p_{1}\right)  +\sum_{i=2}^{k}sgn\left(  p_{i}\right)  g_{i}\right)
=\varepsilon_{\mathbb{Z[\pi}_{1}\left(  B_{2},\ast\right)  ]}\left(
\Gamma^{2},\beta^{n-3}\right)  $. The inclusion of each $g_{i}$ into $\left(
R,\ast\right)  $ yields $j_{\#}\left(  g_{i}\right)  $, thereby converting
$sgn\left(  p_{1}\right)  +\allowbreak\sum_{i=2}^{k}sgn\left(  p_{i}\right)
g_{i}$ to $\widehat{j}_{\#}\left(  sgn\left(  p_{1}\right)  +\allowbreak
\sum_{i=2}^{k}sgn\left(  p_{i}\right)  g_{i}\right)  $. So the lemma is proved.
\end{proof}
\end{lemma}

\section{The Embedded Manifold Plus Construction\label{Sec: EMPC}}

In this section we employ the method of spherical alteration described above
to obtain a constructive proof of Theorem \ref{Th: EPC}. An indirect proof,
relying on the $s$-cobordism theorem, was given in \cite{GT2}. An advantage to
the current approach is that it may be modified to obtain more general
results; that is the content of the last two sections of this paper. A side
benefit of the constructive proof is that it yields a proof of the classical
Manifold Plus Construction which avoids many subtleties related to framings
and bundle theory.

We begin by reviewing a proof of the classical Manifold Plus Construction.
Suppose\emph{ }$\left(  W,A,B\right)  $ is a compact cobordism between closed
manifolds and $A\hookrightarrow W$ is a homotopy equivalence---in which case
we call $\left(  W,A,B\right)  $ is a \emph{one-sided }$h$\emph{-cobordism}.
An application of Poincar\'{e} duality in the universal cover (see
\cite[Th.2.5]{GT1}.) shows that $\pi_{1}\left(  B\right)  \rightarrow\pi
_{1}\left(  W\right)  $ is surjective with perfect kernel. The manifold plus
construction provides a converse to that observation.

\begin{theorem}
[The Manifold Plus Construction]\label{Th: MPC}Let $B$ be a closed $\left(
n-1\right)  $-manifold $\left(  n\geq6\right)  $ and $\theta:\pi_{1}\left(
B,\ast\right)  \rightarrow H$ a surjective homomorphism onto a finitely
presented group such that $\ker\left(  \theta\right)  $ is perfect. Then there
exists a compact one-sided $h$-cobordism $\left(  W,A,B\right)  $ such that
$\ker\left(  \pi_{1}\left(  B,\ast\right)  \rightarrow\pi_{1}\left(
W,\ast\right)  \right)  =\ker\theta$. In fact, it may be arranged that
$A\hookrightarrow W$ is a simple homotopy equivalence, in which case $W$
unique up to homeomorphism rel $B$.
\end{theorem}

A one-sided $h$-cobordism $\left(  W,A,B\right)  $ for which the homotopy
equivalence $A\hookrightarrow W$ is simple will be called a \emph{plus
cobordism}. To avoid repetition, we adopt the convention that whenever a
one-sided [or plus] cobordism is discussed, it will be the first of the two
boundary components listed, i.e., the middle term in the triple, which
includes into $W$ as a [simple] homotopy equivalence.

\begin{proof}
[A classical proof of Theorem \ref{Th: MPC}]\medskip\ \newline\noindent
\textbf{Step I. [Attaching 2-handles to kill }$\ker\left(  \theta\right)
$\textbf{]}\emph{ }Associate $B$ with $B\times\left\{  0\right\}  \subseteq
S=B\times\left[  0,1\right]  $ and let $B_{1}=B\times\left\{  1\right\}  $. By
a standard group theoretic argument, $\ker\theta$ is the normal closure of a
finite set of elements of $\pi_{1}\left(  B,\ast\right)  $; identify a
corresponding collection of nicely embedded oriented loops $\left\{
\alpha_{i}\right\}  _{i=1}^{r}$ in $B_{1}$. Since we are dealing with a normal
closure, we need not be concerned with base points, so we may assume the loops
are pairwise disjoint. Since all elements of $\ker\theta$ are homologically
trivial, each $\alpha_{i}$ has a regular neighborhood in $B_{1}$ homeomorphic
to $S^{1}\times D^{n-2}$. Identify a pairwise disjoint collection $\left\{
J_{i}\right\}  _{i=1}^{r}$ of such neighborhoods and use them as attaching
tubes for a set $\left\{  h_{i}^{2}\right\}  _{i=1}^{r}$ of $2$-handles. For
the moment, we do not concern ourselves with the framings of those
$2$-handles. The resulting $n$-manifold $T$ has fundamental group isomorphic
to $H$. By inverting these handles, $T$ may be obtained by attaching a
collection of $\left(  n-2\right)  $-handles to a collar neighborhood of
$B_{2}=\partial T-B$, a process that does not change fundamental group; so
$B_{2}\hookrightarrow T$ induces a $\pi_{1}$-isomorphism. Note, however, that
$H_{n-2}\left(  T,B_{2};\mathbb{Z}\right)  \cong\mathbb{Z}^{r}$, so we do not
have a one-sided h-cobordism.\medskip

\noindent\textbf{Step II. [Attaching complementary 3-handles]}\emph{ }Here we
will attach a collection of $3$-handles that are complementary to the above
$2$-handles (in an appropriately strong sense) so that the end result is the
desired cobordism $\left(  W,A,B\right)  $. Along the way, we may need to
rechoose the framings of the $2$-handles attached in Step I. Let $\left\{
\beta_{i}^{n-3}\right\}  _{i=1}^{r}$ be the collection of belt spheres of
$\left\{  h_{i}^{2}\right\}  _{i=1}^{r}$. Our initial goal is to identify a
pairwise disjoint collection $\left\{  \Gamma_{i}^{2}\right\}  _{i=1}^{r}$ of
$2$-spheres in $B_{2}$ that is algebraically dual over $\mathbb{Z[\pi}%
_{1}\left(  B_{2},\ast\right)  ]$ to $\left\{  \beta_{i}^{n-3}\right\}
_{i=1}^{r}$. In order to utilize the $\left\{  \Gamma_{i}^{2}\right\}
_{i=1}^{r}$ as attaching spheres for a collection of $3$-handles $\left\{
h_{i}^{3}\right\}  _{i=1}^{r}$, we must also ensure that each has a regular
neighborhood homeomorphic to $S^{2}\times D^{n-3}$ in $B_{2}$. Once that is
accomplished, our task is essentially complete---we will let $W=B\times\left[
0,1\right]  \cup\left\{  h_{i}^{2}\right\}  _{i=1}^{r}\cup\left\{  h_{i}%
^{3}\right\}  _{i=1}^{r}$ and $A=\partial W-B$. Inverting that handle
decomposition, we may view $W$ as a collar on $A$ together with a collection
of $\left(  n-3\right)  $-handles $\left\{  h_{i}^{n-3}\right\}  _{i=1}^{r}$
(the duals of the $3$-handles) with belt spheres $\left\{  \Gamma_{i}%
^{2}\right\}  _{i=1}^{r}$ and a collection of $\left(  n-2\right)  $-handles
$\left\{  h_{i}^{n-2}\right\}  _{i=1}^{r}$ (the duals of the $2$-handles) with
attaching spheres $\left\{  \beta_{i}^{n-3}\right\}  _{i=1}^{r}$. Then
$A\hookrightarrow W$ induces a $\pi_{1}$-isomorphism with each fundamental
group isomorphic to $H$; moreover, the intersection data (as discussed in
Section \ref{subsection: handle theory}) tells us that the corresponding
cellular $\mathbb{Z[\pi}_{1}\left(  A\right)  ]$-complex for the pair $\left(
W,A\right)  $ is of the form
\[
0\rightarrow\widetilde{C}_{n-2}\overset{\partial_{n-2}}{\longrightarrow
}\widetilde{C}_{n-3}\rightarrow0\rightarrow\cdots\rightarrow0
\]
where each of $\widetilde{C}_{n-2}$ and $\widetilde{C}_{n-3}$ is isomorphic to
a free $\mathbb{Z[\pi}_{1}\left(  A\right)  ]$-module on $r$ generators and,
with respect to the obvious preferred bases, the boundary operator
$\partial_{n-2}$ can be represented by a diagonal matrix with diagonal entries
all being $\pm1$. It follows that $A\hookrightarrow W$ is a simple homotopy
equivalence.\smallskip

We now turn to the construction of $\left\{  \Gamma_{i}^{2}\right\}
_{i=1}^{r}$. This is the heart of the matter; it is where the beauty of the
plus construction lies. Since each $\alpha_{i}$ represents an element of the
perfect group $\ker\theta$ it may be expressed as $\prod_{j=1}^{k_{i}}\left[
m_{j}^{i},l_{j}^{i}\right]  $, where $\left\{  \left(  m_{j}^{i},l_{j}%
^{i}\right)  \right\}  _{j=1}^{k_{i}}$ is a complete set of meridian-longitude
pairs for a compact orientable surface $\Lambda_{i}\subseteq B_{1}$ and each
$m_{j}^{i}$ and $l_{j}^{i}$ also lies in $\ker\theta$. Using general position
and the radial structure of the attaching tubes $J_{i}$, we can adjust these
surfaces so that each has boundary $\alpha_{i}^{\prime}$ which lies in
$\partial J_{i}$ and is parallel to $\alpha_{i}$. In addition we may assume
that the $\Lambda_{i}$ are properly embedded in $B_{1}-%
{\textstyle\bigcup\nolimits_{i=1}^{r}}
\operatorname{int}J_{i}$ and pairwise disjoint. Complete each $\Lambda_{i}$ to
a closed surface $\widehat{\Lambda}_{i}\subseteq B_{2}$ by adding a $2$-disk
$\widehat{D}_{i}$ lying in the belt tube of $h_{i}^{2}$ and parallel to its
core. It is here that we must pay attention to framings. Since $\Lambda_{i}$
is orientable and deformation retracts onto a bouquet of circles, where each
of those circles corresponds to an $m_{j}^{i}$ or an $l_{j}^{i}$---all of
which have trivial normal bundles, standard bundle theory can be used to
verify that $\Lambda_{i}$ has a product relative regular neighborhood in
$B_{1}-%
{\textstyle\bigcup\nolimits_{i=1}^{r}}
\operatorname{int}J_{i}$. If necessary, we now rechoose the framing used to
attach $h_{i}^{2}$ so that the corresponding trivial normal bundle for
$\widehat{D}_{i}$ matches up with that of $\Lambda_{i}$ to give $\widehat
{\Lambda}_{i}$ a product regular neighborhood in $B_{2}$; indeed, this is
precisely the matter which determines the framings that must be used for
attaching the $2$-handles. (One may argue that this should have been discussed
\emph{before} attaching the $2$-handles; however, it seems instructive to
discover the issue in context.)

Notice that each surface $\widehat{\Lambda}_{i}$ intersects the belt sphere
$\beta_{i}^{n-3}$ transversely in exactly one point and that it intersects no
other belt spheres. Thus, $\left\{  \widehat{\Lambda}_{i}\right\}  $ is a
collection of geometric duals\ for $\left\{  \beta_{i}^{n-3}\right\}  $, and
since the fundamental group of each $\widehat{\Lambda}_{i}$ includes trivially
into $B_{2}$, this may be expressed in terms of $\mathbb{Z[\pi}_{1}\left(
B_{2}\right)  ]$-intersection numbers. After choosing all necessary base
points, base paths, and orientations the geometric intersection properties
imply that $\varepsilon_{\mathbb{Z[\pi}_{1}\left(  B_{2}\right)  ]}\left(
\widehat{\Lambda}_{i},\beta_{j}^{n-3}\right)  =0$ whenever $i\neq j$ and each
$\varepsilon_{\mathbb{Z[\pi}_{1}\left(  B_{2}\right)  ]}\left(  \widehat
{\Lambda}_{i},\beta_{i}^{n-3}\right)  =\pm g_{i}$ for some $g_{i}%
\in\mathbb{Z[\pi}_{1}\left(  B_{2}\right)  ]$. We may arrange that each of the
latter intersection numbers is $\pm1$ by rechoosing some of the base paths.

Unfortunately, the $\widehat{\Lambda}_{i}$ will usually have genus $>0$, and
thus be unusable for attaching $3$-handles. We remedy that problem by
surgering the surfaces to $2$-spheres in the manner outlined in Section
\ref{subsection: surgery on surfaces}. Since each $m_{j}^{i}$ contracts in
$B_{2}$ we may surger $\widehat{\Lambda}_{i}$ to a $2$-sphere $\Gamma_{i}$ in
$B_{2}$ using disks bounded by the various meridianal curves. (In the notation
of Section \ref{subsection: surgery on surfaces}, $\Gamma_{i}=\widehat
{\Lambda}_{i}^{\ast}$.) By part 1) of Lemma \ref{Lemma: Surgery on surfaces}
this operation preserves $\mathbb{Z[\pi}_{1}\left(  B_{2}\right)
]$-intersection numbers, so $\varepsilon_{\mathbb{Z[\pi}_{1}\left(
B_{2}\right)  ]}\left(  \beta_{i}^{n-3},\Gamma_{j}\right)  =\pm\delta_{ij}$
for all $1\leq i,j\leq r$. Another application of standard bundle theory
ensures that the $\Gamma_{i}$ inherit trivial normal bundles from the
$\widehat{\Lambda}_{i}$, so they may be used as attaching spheres for the
$3$-handles $\left\{  h_{i}^{3}\right\}  _{i=1}^{r}$, thereby supplying the
final ingredient of the manifold plus construction.

The uniqueness part of this theorem follows from a clever application of the
$s$-cobordism theorem. Since it is not of primary importance to this paper, we
refer the reader to \cite[p.197]{FQ} for a proof.
\end{proof}

We are now ready for the embedded version of the Manifold Plus Construction.
Much of the strategy and notation employed above is recycled into the
proof---the main ideas are the same. Some issues become more complex due to
our desire to embed the construction in an ambient manifold; as a pleasant
surprise, other issues become easier for the same reason.

\begin{theorem}
[The Embedded Manifold Plus Construction]\label{Th: EPC}Let $R$ be an
$n$-manifold ($n\geq6$) containing a closed $\left(  n-1\right)  $-manifold
$B$ in its boundary and suppose \linebreak$\ker\left(  \pi_{1}\left(
B,\ast\right)  \overset{i_{\ast}}{\longrightarrow}\pi_{1}\left(
R,\ast\right)  \right)  $ contains a perfect group $G$ which is the normal
closure in $\pi_{1}\left(  R,\ast\right)  $ of a finite set of elements. Then
there exists an embedding of a plus cobordism $\left(  W,A,B\right)  $ into
$R$ which is the identity on $B$ and for which $\ker\left(  \pi_{1}\left(
B\right)  \rightarrow\pi_{1}\left(  W\right)  \right)  =G$.

\begin{proof}
\noindent\textbf{Step I.}\emph{ }\textbf{ [Finding embedded 2-handles that
kill }$\ker\left(  i_{\#}\right)  $\textbf{] \ }Let $S\approx B\times\left[
0,1\right]  $ be a collar neighborhood of $B$ in $R$ and let $B_{1}$ denote
the interior boundary component of $S$. Choose a pairwise disjoint collection
of properly embedded $2$-disks $\left\{  D_{1},\cdots,D_{r}\right\}  $ in
$\overline{R-S}$ whose boundaries in $B_{1}$ represent a finite normal
generating set for $G$. By taking regular neighborhoods, thicken the $D_{i}$
to a pairwise disjoint collection of $2$-handles $\{h_{i}^{2}\}_{i=1}^{r}$.
Let $T=S\cup$ $h_{1}^{2}\cup\cdots\cup h_{r}^{2}$ and $B_{2}=\partial T-B$;
for later use, let $J_{i}$ denote the attaching tube for $h_{i}^{2}$. Then
$\pi_{1}\left(  B_{2}\right)  \cong\pi_{1}\left(  T\right)  \cong\pi
_{1}\left(  B\right)  /G$ and $\ker\left(  \pi_{1}\left(  B\right)
\rightarrow\pi_{1}\left(  T\right)  \right)  =G$. For the remainder of the
proof, all work will be done within a regular neighborhood $R^{\prime}$ of $T$
in $R$. Since $R^{\prime}$ is just $T$ with a collar added along $B_{2}$,
$B_{2}\hookrightarrow\overline{R^{\prime}-T}$ induces a $\pi_{1}%
$-isomorphism---a fact that will be utilized only in the special argument
needed for the $n=6$ case.\medskip

\noindent\textbf{Step II.}\emph{ }\textbf{[Altering the embedded 2-handles so
that complementary embedded 3-handles exist] \ }We would like to find a
pairwise disjoint collection of $3$-handles $\left\{  h_{i}^{3}\right\}
_{i=1}^{r}$ embedded $\overline{R^{\prime}-T}$ with attaching $2$-spheres
algebraically dual over $\mathbb{Z[\pi}_{1}\left(  B_{2}\right)  ]$ to the
collection $\left\{  \beta_{i}^{n-3}\right\}  _{i=1}^{r}$of belt spheres of
$\left\{  h_{i}^{2}\right\}  _{i=1}^{r}$.\emph{ }Adding those handles to $T$
(and following the argument used in the previous theorem) would give us the
desired $W\subseteq R^{\prime}$.

Toward that end goal, we construct a collection $\left\{  \Gamma_{i}\right\}
_{i=1}^{r}$ of $2$-spheres in $B_{2}$ which are algebraic duals for the
collection $\left\{  \beta_{i}^{n-3}\right\}  _{i=1}^{r}$ in precisely the
same manner employed in Step II of the previous theorem. (But unlike that
proof, we need not concern ourselves with framings of the $2$-handles or
regular neighborhoods of the $2$-spheres.)

Under ideal circumstances, the $\left\{  \Gamma_{i}\right\}  _{i=1}^{r}$ would
contract in $\overline{R^{\prime}-T}$ allowing us to obtain a pairwise
disjoint collection of properly embedded $3$-disks in $\overline{R^{\prime}%
-T}$ with the $\Gamma_{i}$ as boundaries. Thickening those disks to
$3$-handles would complete the construction of $W$.

The main strategy of this proof can now be described: by utilizing a carefully
selected sequence of spherical modifications of the $\left\{  h_{i}%
^{2}\right\}  _{i=1}^{r}$ we arrive at a new collection of embedded
$2$-handles so that the correspondingly altered versions of the $\left\{
\Gamma_{i}\right\}  _{i=1}^{r}$ satisfy the desired contractibility
condition.\medskip

\noindent\emph{Step (II}$_{1}$\emph{)}\textsc{. }(Spherical alteration of
$h_{1}^{2}.$) Locate a parallel copy\ of $\Gamma_{1}$ lying in the interior of
$\overline{R^{\prime}-T}$. (Push $\Gamma_{1}$ along collar lines.) Reverse the
orientation on that copy and denote it by $\Delta_{1}$. Choose a base point
$p_{1}$ for $\Gamma_{1}$ lying in the belt tube of $h_{1}^{2}$ but missing the
belt sphere $\beta_{1}^{n-3}$, and let $\xi_{1}$ be the track of $p_{1}$ in
$\overline{R^{\prime}-T}$ under the push. Extend $\xi_{1}$ slightly to an arc
$\xi_{1}^{\prime}$ which connects the core of $h_{1}^{2}$ to $\Delta_{1}$, and
perform a spherical alteration of $h_{1}^{2}$ over $\Delta_{1}$ along $\xi
_{1}^{\prime}$ to obtain $h_{1}^{2}\left(  \xi_{1}^{\prime},\Delta_{1}\right)
$. Then perform the corresponding spherical alterations on each $2$-sphere in
the collection $\left\{  \Gamma_{i}\right\}  _{i=1}^{r}$ (\textbf{Note.
}Although the \emph{algebraic} intersection number of $\Gamma_{i}$ with
$\beta_{1}^{n-3}$ is $0$ when $i\geq2$, the two need not be disjoint; so the
alterations must be done in order to obtain a collection that lies in the
right-hand boundary $B_{2}^{\left(  1\right)  }$ of $T^{\left(  1\right)
}=S\cup h_{1}^{2}\left(  \xi_{1}^{\prime},\Delta_{1}\right)  \cup h_{2}%
^{2}\cup\cdots h_{r}^{2}$.) Now make the following observations:

\begin{enumerate}
\item[a$_{1}$)] The collection $\left\{  \Gamma_{i}^{2}\left(  \xi_{1}%
^{\prime},\Delta_{1}\right)  \right\}  _{i=1}^{r}$ of altered $2$-spheres is
algebraically dual to the set $\{\beta_{i}^{n-3}\}_{i=1}^{r}$ of belt spheres
of the new collection of $2$-handles $\{h_{1}^{2}\left(  \xi_{1}^{\prime
},\Delta_{1}\right)  ,\allowbreak h_{2}^{2},\allowbreak\cdots,h_{r}^{2}\}$ in
$B_{2}^{\left(  1\right)  }$. (Recall that the belt sphere $\beta_{1}^{n-3}$
of $h_{1}^{2}$ is also the belt sphere for $h_{1}^{2}\left(  \xi_{1}^{\prime
},\Delta_{1}\right)  $.) We need only check that $\varepsilon_{\mathbb{Z[\pi
}_{1}\left(  B_{2}^{\left(  1\right)  },\ast\right)  ]}\left(  \beta_{i}%
^{n-3},\Gamma_{j}^{2}\left(  \xi_{1}^{\prime},\Delta_{1}\right)  \right)
=\varepsilon_{\mathbb{Z\mathbb{[}\pi}_{1}\left(  B_{2},\ast\right)  ]}\left(
\beta_{i}^{n-3},\Gamma_{j}\right)  $ for all $1\leq i,j\leq r$. But this is
clear, since $\Gamma_{i}^{2}\left(  \xi_{1}^{\prime},\Delta_{1}\right)  $ and
$\Gamma_{i}$ are identical over path connected subsets which contain all
points of intersection with elements of $\{\beta_{i}^{n-3}\}_{i=1}^{r}$. This
means that all paths and loops utilized in determining the two intersection
numbers can be chosen to be identical; such loops represent the
`same'\ elements of $\pi_{1}\left(  B_{2}^{\left(  1\right)  }\right)  $ as
they do in $\pi_{1}\left(  B_{2}\right)  $.

\item[b$_{1}$)] $\Gamma_{1}^{2}\left(  \xi_{1}^{\prime},\Delta_{1}\right)  $
contracts in $R^{\prime}$. To see this, let $\lambda_{1}$ be a path from the
base point $p_{0}$ of $R^{\prime}$ to $p_{1}$, then by Lemma
\ref{alteration lemma}, as elements of $\pi_{2}\left(  R^{\prime}%
,p_{0}\right)  $,
\begin{align*}
\left[  \lambda_{1}\Gamma_{1}^{2}\left(  \xi_{1}^{\prime},\Delta_{1}\right)
\right]   &  =\left[  \lambda_{1}\Gamma_{1}^{2}\right]  +\left[  \left(
\lambda_{1}\cdot\xi_{1}\right)  \left(  \Delta_{1}\right)  \right] \\
&  =\left[  \lambda_{1}\Gamma_{1}^{2}\right]  +\left[  \lambda_{1}\left(
-\Gamma_{1}^{2}\right)  \right] \\
&  =\left[  \lambda_{1}\Gamma_{1}^{2}\right]  -\left[  \lambda_{1}\left(
\Gamma_{1}^{2}\right)  \right]  =0.
\end{align*}
So $\Gamma_{1}^{2}\left(  \xi_{1}^{\prime},\Delta_{1}\right)  $ contracts in
$R^{\prime}$.\medskip
\end{enumerate}

\noindent\emph{Step (II}$_{2}$\emph{)}\textsc{. }(Spherical alteration of
$h_{2}^{2}.$) Now begin with the collection of $2$-handles $\left\{  h_{1}%
^{2}\left(  \xi_{1}^{\prime},\Delta_{1}\right)  ,h_{2}^{2},\cdots,h_{r}%
^{2}\right\}  $ attached to $S$ in $R^{\prime}$ and the collection $\left\{
\Gamma_{i}^{2}\left(  \xi_{1}^{\prime},\Delta_{1}\right)  \right\}  _{i=1}%
^{r}$ of algebraic duals to their belt spheres $\{\beta_{i}^{n-3}\}_{i=1}^{r}$
in $B_{2}^{\left(  1\right)  }$. Obtain a $2$-sphere $\Delta_{2}$ by pushing
$\Gamma_{2}\left(  \xi_{1}^{\prime},\Delta_{1}\right)  $ into the interior of
$\overline{R^{\prime}-T^{\left(  1\right)  }}$ and then reversing its
orientation. Let $p_{2}$ be an appropriately chosen base point for $\Gamma
_{2}\left(  \xi_{1}^{\prime},\Delta_{1}\right)  $ and $\xi_{2}$ be the track
of $p_{2}$ in $\overline{R^{\prime}-T^{\left(  1\right)  }}$ under the push.
Extend $\xi_{2}$ slightly to an arc $\xi_{2}^{\prime}$ which connects the core
of $h_{2}^{2}$ to $\Delta_{2}$, and perform the spherical alteration of
$h_{2}^{2}$ over $\Delta_{2}$ along $\xi_{2}^{\prime}$ to obtain $h_{2}%
^{2}\left(  \xi_{2}^{\prime},\Delta_{2}\right)  $. Then perform the
corresponding spherical alterations on each $2$-sphere in the collection
$\left\{  \Gamma_{i}^{2}\left(  \xi_{1}^{\prime},\Delta_{1}\right)  \right\}
_{i=1}^{r}$. To save on notation, we denote the twice altered version of
$\Gamma_{i}^{2}$ by $\Gamma_{i}^{2}\left(  \sqcup_{j=1}^{2}\xi_{j}^{\prime
},\Delta_{j}\right)  $. Let $T^{\left(  2\right)  }=S\cup h_{1}^{2}\left(
\xi_{1}^{\prime},\Delta_{1}\right)  \cup h_{2}^{2}\left(  \xi_{2}^{\prime
},\Delta_{2}\right)  \cup h_{3}^{2}\cdots h_{r}^{2}$. Using the same arguments
as above, we now have:

\begin{enumerate}
\item[a$_{2}$)] The collection $\left\{  \Gamma_{i}^{2}\left(  \sqcup
_{j=1}^{2}\xi_{j}^{\prime},\Delta_{j}\right)  \right\}  _{i=1}^{r}$ of twice
altered $2$-spheres is algebraically dual to the set $\{\beta_{i}%
^{n-3}\}_{i=1}^{r}$ of belt spheres of the new collection of $2$-handles
$\{h_{1}^{2}\left(  \xi_{1}^{\prime},\Delta_{1}\right)  ,h_{2}^{2}\left(
\xi_{2}^{\prime},\Delta_{2}\right)  ,h_{3}^{2},\cdots,h_{r}^{2}\}$ in
$B_{2}^{\left(  2\right)  }$.

\item[b$_{2}$)] $\Gamma_{1}^{2}\left(  \sqcup_{j=1}^{2}\xi_{j}^{\prime}%
,\Delta_{j}\right)  $ and $\Gamma_{2}^{2}\left(  \sqcup_{j=1}^{2}\xi
_{j}^{\prime},\Delta_{j}\right)  $ both contract in $R^{\prime}$.
(Contractibility of $\Gamma_{1}^{2}\left(  \sqcup_{j=1}^{2}\xi_{j}^{\prime
},\Delta_{j}\right)  $ follows from Lemma \ref{alteration lemma} together with
the fact that \newline$\varepsilon_{\mathbb{Z\pi}_{1}\left(  B_{2}^{\left(
1\right)  }\right)  }\left(  \beta_{2}^{n-3},\allowbreak\Gamma_{1}^{2}\left(
\xi_{1}^{\prime},\Delta_{1}\right)  \right)  =0$.)
\end{enumerate}

Continue the above process until each of $h_{1}^{2},\cdots,h_{r}^{2}$ has been
altered by a similar process, and let $T^{\left(  r\right)  }$ denote the
union of $S$ and these altered $2$-handles. At this point we have:

\begin{enumerate}
\item[a$_{r}$)] The collection $\left\{  \Gamma_{i}^{2}\left(  \sqcup
_{j=1}^{r}\xi_{j}^{\prime},\Delta_{j}\right)  \right\}  _{i=1}^{r}$ of
$r$-times altered $2$-spheres is algebraically dual to the set $\{\beta
_{i}^{n-3}\}_{i=1}^{r}$of belt spheres of the handles $\{h_{i}^{2}\left(
\xi_{i}^{\prime},\Delta_{i}\right)  \}_{i=1}^{r}$ in $B_{2}^{\left(  r\right)
}$.

\item[b$_{r}$)] $\Gamma_{1}^{2}\left(  \sqcup_{j=1}^{r}\xi_{j}^{\prime}%
,\Delta_{j}\right)  ,\cdots,\Gamma_{r}^{2}\left(  \sqcup_{j=1}^{r}\xi
_{j}^{\prime},\Delta_{j}\right)  $ contract in $R^{\prime}$.
\end{enumerate}

Before proceeding to the final stage of our proof, let us simplify the above
notation. From now on each of the altered $2$-handles $h_{i}^{2}\left(
\xi_{i}^{\prime},\Delta_{i}\right)  $ will be denoted by $\dot{h}_{i}^{2}$ and
each of the $r$-times altered $2$-spheres $\Gamma_{i}^{2}\left(  \sqcup
_{j=1}^{r}\xi_{j}^{\prime},\Delta_{j}\right)  $ by $\dot{\Gamma}_{i}^{2}$ .
Thus, we have $T^{\left(  r\right)  }=S\cup\dot{h}_{1}^{2}\cup\cdots\cup
\dot{h}_{r}^{2}\subseteq R^{\prime}$ and a collection $\left\{  \dot{\Gamma
}_{i}^{2}\right\}  _{i=1}^{r}$of $2$-spheres algebraically dual over
$\mathbb{Z[\pi}_{1}\left(  B_{2}^{\left(  r\right)  },\ast\right)  ]$ to the
collection $\left\{  \beta_{i}^{n-3}\right\}  _{i=1}^{r}$ of belt spheres of
those $2$-handles in $B_{2}^{\left(  r\right)  }$. In addition, each
$\dot{\Gamma}_{i}^{2}$ contracts in $R^{\prime}$. By general position, these
$2$-spheres also contract in $\overline{R^{\prime}-T^{\left(  r\right)  }}$.
This is because any contraction of $\dot{\Gamma}_{i}^{2}$ in $R^{\prime}$ can
be pushed off the $2$-dimensional cores of all of the $2$-handles and, thus,
entirely out of the interior of $T^{\left(  r\right)  }$. Assume for the
moment that the dimension of $R^{\prime}$ is at least $7$, Then we may choose
a pairwise disjoint collection $\left\{  D_{i}^{3}\right\}  _{i=1}^{r}$ of
properly embedded $3$-disks in $\overline{R^{\prime}-T^{\left(  r\right)  }}$
such that $\partial D_{i}^{3}=\dot{\Gamma}_{i}^{2}$ for each $i=1,\cdots,r$.
Take pairwise disjoint relative regular neighborhoods of these $3$-disks in
$\overline{R^{\prime}-T^{\left(  r\right)  }}$ to obtain $3$-handles $\left\{
\dot{h}_{i}^{3}\right\}  _{i=1}^{r}$ in $\overline{R^{\prime}-T^{\left(
r\right)  }}$ attached to $B_{2}^{\left(  r\right)  }$. Let
\[
W=T^{\left(  r\right)  }\cup\left(  \cup_{i=1}^{r}\dot{h}_{i}^{3}\right)
=S\cup\left(  \cup_{i=1}^{r}\dot{h}_{i}^{2}\right)  \cup\left(  \cup_{i=1}%
^{r}\dot{h}_{i}^{3}\right)
\]
and let $A=\partial W-B$. By the same reasoning used in Theorem \ref{Th: MPC},
$\left(  W,A,B\right)  $ is a plus cobordism.\medskip

\textbf{Step III. [}$\pi$-$\pi$ \textbf{argument for the }$n=6$\textbf{
case.]} The only place the above proof runs into trouble is in the use of
general position to obtain a pairwise disjoint collection $\left\{  D_{i}%
^{3}\right\}  _{i=1}^{r}$ of properly embedded $3$-disks in $\overline
{R^{\prime}-T^{\left(  r\right)  }}$ with $\partial D_{i}^{3}=\dot{\Gamma}%
_{i}^{2}$ for each $i$. If $n=6$, we may use general position to obtain a
collection $\left\{  \widetilde{D}_{i}^{3}\right\}  _{i=1}^{r}$of
\emph{immersed} $3$-disks, each containing a finite collection of interior
transverse self-intersection points, and a finite number of interior points
where it transversely intersects another member of the collection. We will
employ a well known strategy (see, for example, the proof of the $\pi$-$\pi$
Theorem in \cite[Ch.4]{Wa}) to eliminate all intersection and
self-intersection points. Once that is accomplished, the proof may be
completed in the previous manner.

Let $p\in\widetilde{D}_{i}^{3}\cap\widetilde{D}_{j}^{3}$ for some $i\neq j$.
Choose arcs $\alpha$ and $\alpha^{\prime}$ in $\widetilde{D}_{i}^{3}$ and
$\widetilde{D}_{j}^{3}$ respectively missing all other intersection points and
connecting $p$ to points in $B_{2}^{\left(  r\right)  }$. Since $B_{2}%
^{\left(  r\right)  }\hookrightarrow\overline{R^{\prime}-T^{\left(  r\right)
}}$ induces a $\pi_{1}$-isomorphism (a surjection is sufficient), we may
connect the endpoints of $\alpha$ and $a^{\prime}$ by an arc $\alpha
^{\prime\prime}$ in $B_{2}^{\left(  r\right)  }$ such that the loop
$\alpha\cup a^{\prime}\cup\alpha^{\prime\prime}$ contracts in $\overline
{R^{\prime}-T^{\left(  r\right)  }}$. By general position we may choose an
embedded $2$-disk $\delta$ bounded by $\alpha\cup a^{\prime}\cup\alpha
^{\prime\prime}$ which intersects the collection $\left\{  \widetilde{D}%
_{k}^{3}\right\}  _{k=1}^{r}$ only at $\alpha\cup\alpha^{\prime}$; and
intersects $B_{2}^{\left(  r\right)  }$ only at $\alpha^{\prime\prime}$. Use
$\delta$ to define a proper isotopy of $\widetilde{D}_{i}^{3}$ in
$\overline{R^{\prime}-T^{\left(  r\right)  }}$ which moves points in a small
neighborhood of $\alpha$ across $\delta$ to the other side of $\alpha^{\prime
}$, thus eliminating the point of intersection $p$. See Figure \ref{Fig4}.%
\begin{figure}
[ptb]
\begin{center}
\includegraphics[
height=2.7332in,
width=5.5536in
]%
{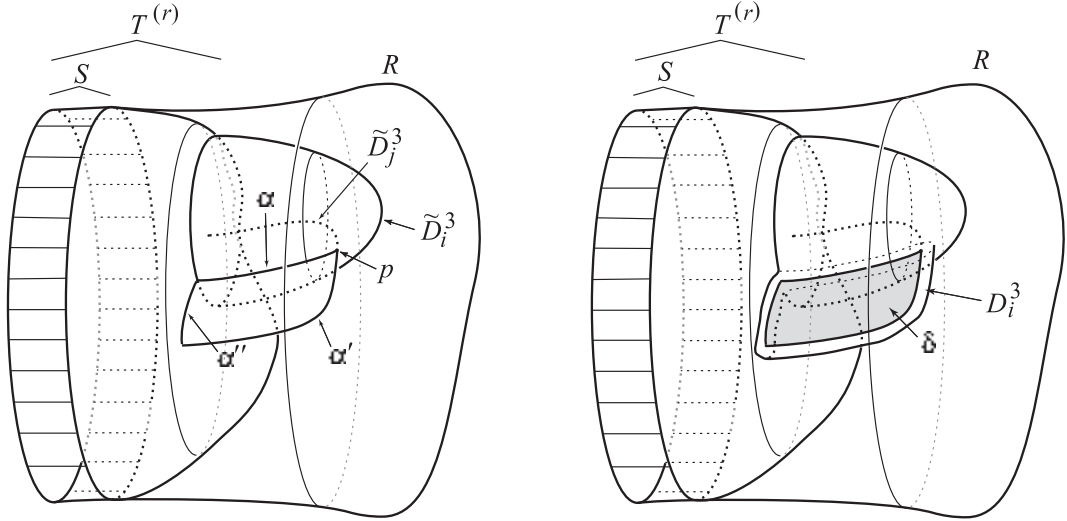}%
\caption{Eliminating a point of intersection using the $\pi$-$\pi$ procedure}%
\label{Fig4}%
\end{center}
\end{figure}
Alternatively, perform a \textquotedblleft finger move\textquotedblright\ on
$\widetilde{D}_{i}^{3}$ along the arc $\alpha^{\prime\prime}$ to create a new
point $p^{\prime}$ of transverse intersection between $\widetilde{D}_{i}^{3}$
and $\widetilde{D}_{j}^{3}$. The disk $\delta$ contains a Whitney disk which
allows us to simultaneously remove $p$ and $p^{\prime}$ via an isotopy of
$\widetilde{D}_{i}^{3}$ (after finger move) which takes place entirely in a
neighborhood of $\delta$. Similar procedures may be used to eliminate a point
of self-intersection from a given $\widetilde{D}_{i}^{3}$.

Apply the above procedure to each point of intersection between distinct
elements of $\left\{  \widetilde{D}_{k}^{3}\right\}  _{k=1}^{r}$ and each
point of self-intersection of each $\widetilde{D}_{i}^{3}$ to arrive at a
pairwise disjoint collection of properly embedded $3$-disks $\left\{
D_{k}^{3}\right\}  _{k=1}^{r}$. Note that the boundary of each $D_{k}^{3}$ is
isotopic to the original $\dot{\Gamma}_{k}^{2}$ in $B_{2}^{\left(  r\right)
}$. Hence, the $\left\{  \partial D_{k}^{3}\right\}  _{k=1}^{r}$ is still a
collection of algebraic duals for the $\left\{  \beta_{k}^{n-3}\right\}
_{k=1}^{r}$ in $B_{2}^{\left(  r\right)  }$.
\end{proof}
\end{theorem}

\begin{remark}
\label{Remark: EPC implies MPC}\emph{The reader will note that a key issue in
the proof of Theorem \ref{Th: MPC}---the existence of product neighborhoods
for the }$2$\emph{-spheres along which the }$3$\emph{-handles will be
attached---does not appear in the proof of Theorem \ref{Th: EPC}. In the
latter setting, the }$3$\emph{-handles are realized as regular neighborhoods
of embedded }$3$\emph{-disks; as such, product neighborhoods of their
boundaries are guaranteed by regular neighborhood theory.}

\emph{This is the essence of our alternate proof of Theorem \ref{Th: MPC}. One
first carries out Step I of the classical proof; in particular, construct a
manifold }$T$\emph{ by attaching finitely many }$2$\emph{-handles to }%
$B\times\left[  0,1\right]  $\emph{ to kill }$\ker\theta$\emph{ (and with no
attention given to the framings used). Theorem \ref{Th: EPC} applied to the
inclusion }$B\overset{i}{\hookrightarrow}T$\emph{ with }$G=\ker i_{\ast}%
=\ker\theta$\emph{ then assures the existence of the desired plus cobordism
lying inside }$T$\emph{. It strikes us as surprising that the full plus
cobordism can be found embedded in the first stage of that construction---even
when the first stage is done with the wrong framings.}
\end{remark}

\section{Generalized manifold plus constructions and their embeddings}

The techniques employed in the proofs of Theorems \ref{Th: MPC} and
\ref{Th: EPC} can be carried out without the full hypothesis of `perfectness'
on the subgroups $\ker\theta$ and $G$, provided one is satisfied with weaker
(but still useful) conclusions. In this section we develop results of that
type. Primary motivation for the definitions and results found here is
provided by our ongoing study of noncompact manifolds \cite{GT3}.

Our first goal is to formulate appropriate generalizations of `one-sided
h-cobordism' and `plus cobordism'. Let $\left(  X,A\right)  $ be a connected
CW pair and $L\trianglelefteq\pi_{1}\left(  A\right)  $. The inclusion
$A\hookrightarrow X$ is a ($\operatorname{mod}L$)\emph{-homotopy equivalence}
if it induces an isomorphism on fundamental groups and is a $\mathbb{Z}\left[
\pi_{1}\left(  A\right)  /L\right]  $-homology equivalence; if that homology
equivalence is simple we call $A\hookrightarrow X$ is a ($\operatorname{mod}%
L$)\emph{-simple homotopy equivalence}.

A compact cobordism $\left(  W,A,B\right)  $ is a ($\operatorname{mod}%
L$)\emph{-one-sided h-cobordism} if $B\hookrightarrow W$ induces a surjection
of fundamental groups and $A\hookrightarrow W$ is a ($\operatorname{mod}%
L$)-homotopy equivalence. A one-sided ($\operatorname{mod}L$)-h-cobordism for
which $A\hookrightarrow W$ is a ($\operatorname{mod}L$)-simple homotopy
equivalence is called a ($\operatorname{mod}L$)\emph{-plus cobordism}.

\begin{remark}
\emph{Standard arguments show that the notions of (}$\operatorname{mod}%
L$\emph{)-homotopy equivalence, (}$\operatorname{mod}L$\emph{)-one-sided
h-cobordism, and (}$\operatorname{mod}L$\emph{)-plus cobordism reduce to the
classical definitions when }$L=\left\{  1\right\}  $\emph{. In that case the
surjectivity of }$\pi_{1}\left(  B\right)  \rightarrow\pi_{1}\left(  W\right)
$\emph{ is automatic and need not be included in the definitions (see the
second paragraph of Section \ref{Sec: EMPC}). When }$L$\emph{ is nontrivial
that condition must be included in the definition to obtain a theory that
parallels the classical situation. For example, it provides a natural
correspondence between }$\pi_{1}\left(  A\right)  /L$\emph{ and }$\pi
_{1}\left(  B\right)  /L^{\prime}$\emph{ where }$L^{\prime}$\emph{ is the
preimage of }$L$\emph{; from there it follows (by Poincar\'{e} duality) that
}$B\hookrightarrow W$\emph{ is also a }$Z\left[  \pi_{1}\left(  A\right)
/L\right]  $\emph{-homology equivalence.}
\end{remark}

The following provides an important connection between ($\operatorname{mod}%
L$)-one-sided h-cobordisms and the material presented in Section
\ref{Subsection: perfect and nearly perfect groups}.

\begin{lemma}
Let $\left(  W,A,B\right)  $ is a ($\operatorname{mod}L$)-one-sided
h-cobordism and $i_{\#}:\pi_{1}\left(  B\right)  \rightarrow\pi_{1}\left(
W\right)  $ the inclusion induced surjection. Then $\ker i_{\#}$ is strongly
$L^{\prime}$-perfect, where $L^{\prime}=i_{\#}^{-1}\left(  L\right)  $.

\begin{proof}
This follows from Poincar\'{e} duality and the 5-Term Exact Sequence from the
theory of group homology \cite{Stal}, \cite{Stam}. See \cite{GT3} for a
detailed proof.
\end{proof}
\end{lemma}

We are now ready to state and prove generalizations of the two main theorems
from the previous section.

\begin{theorem}
[A Generalized Manifold Plus Construction]\label{Th: GMPC}Let $B$ be a
closed\linebreak\ $\left(  n-1\right)  $-manifold $\left(  n\geq6\right)  $
and $\theta:\pi_{1}\left(  B,\ast\right)  \rightarrow H$ a surjective
homomorphism onto a finitely presented group such that $\ker\left(
\theta\right)  $ is strongly $L^{\prime}$-perfect for some group $L^{\prime}$
where $\ker\left(  \theta\right)  \trianglelefteq L^{\prime}\trianglelefteq
\pi_{1}\left(  B,\ast\right)  $ and $\omega_{1}\left(  L^{\prime}\right)  =1$.
Then, for $L=L^{\prime}/\ker\theta$, there exists a compact
($\operatorname{mod}L$)-plus cobordism $\left(  W,A,B\right)  $ such that
$\ker\left(  \pi_{1}\left(  B,\ast\right)  \rightarrow\pi_{1}\left(
W,\ast\right)  \right)  =\ker\theta$.
\end{theorem}

\begin{theorem}
[A Generalized Embedded Manifold Plus Construction]\label{Th: GEPC}Let $R$ be
an $n$-manifold ($n\geq6$) containing a closed $\left(  n-1\right)  $-manifold
$B$ in its boundary and suppose $\pi_{1}\left(  B,\ast\right)  $ contains a
pair of normal subgroups $G\leq L^{\prime}$, each contained in $\ker\left(
\pi_{1}\left(  B,\ast\right)  \overset{i_{\ast}}{\longrightarrow}\pi
_{1}\left(  R,\ast\right)  \right)  $, such that $G$ is strongly $L^{\prime}%
$-perfect. Suppose also that $G$ is the normal closure in $\pi_{1}\left(
B,\ast\right)  $ of a finite set of elements. Then, for $L=L^{\prime}/G$,
there exists an embedding of a ($\operatorname{mod}L$)-plus cobordism $\left(
W,A,B\right)  $ into $R$ which is the identity on $B$ and for which
$\ker\left(  \pi_{1}\left(  B\right)  \rightarrow\pi_{1}\left(  W\right)
\right)  =G$.
\end{theorem}

Proofs of each of these theorems can be obtained by reworking those from the
previous section with the new weaker hypotheses---obtaining correspondingly
weaker conclusions. We sketch out the details of those changes needed to
obtain Theorem \ref{Th: GMPC} and leave it to the reader to carry out the
analogous changes required to obtain Theorem \ref{Th: GEPC}.

\begin{proof}
[Sketch of the Generalized Manifold Plus Construction]Begin by repeating Step
I of the proof of Theorem \ref{Th: MPC}---in particular, attach a collection
$\left\{  h_{i}^{2}\right\}  _{i=1}^{r}$of $2$-handles to $B_{1}%
=B\times\left\{  1\right\}  \subseteq S=B\times\left[  0,1\right]  $ to kill a
finite set $\left\{  \alpha_{i}\right\}  _{i=1}^{r}$ of oriented loops which
normally generate $\ker\left(  \theta\right)  $. Let $T=S\cup(\cup_{i=1}%
^{r}h_{i}^{2})$. Moving to Step II, strong $L^{\prime}$-perfectness ensures
that each $\alpha_{i}$ bounds a compact oriented surface $\Lambda_{i}\subseteq
B_{1}$ which contains a complete set $\left\{  \left(  m_{j}^{i},l_{j}%
^{i}\right)  \right\}  _{j=1}^{k_{i}}$ of meridian-longitude pairs for which
each $m_{j}^{i}$ corresponds to an element of $\ker\theta$ and each $l_{j}%
^{i}$ to an element of $L^{\prime}$. Using the same general position and
bundle-theoretic arguments employed earlier, add a disk to each $\Lambda_{i}$
to obtain a pairwise disjoint collection $\left\{  \widehat{\Lambda}%
_{i}\right\}  _{i=1}^{r}$ of closed oriented surfaces in $B_{2}=\partial T-B$
geometrically dual to the collection $\left\{  \beta_{i}^{n-3}\right\}
_{i=1}^{r}$ of belt spheres of the $\left\{  h_{i}^{2}\right\}  _{i=1}^{r}$.
By rechoosing the framings of the $2$-handles if necessary, we may arrange
that each $\widehat{\Lambda}_{i}$ has a product neighborhood in $B_{2}$---here
we utilize the hypothesis that $\omega\left(  L^{\prime}\right)  =1$. Since
each $\widehat{\Lambda}_{i}$ has a fundamental group that includes into
$L\trianglelefteq\mathbb{\pi}_{1}\left(  B_{2}\right)  $, the $\mathbb{Z[\pi
}_{1}\left(  B_{2}\right)  /L]$-intersection numbers between elements of
$\left\{  \widehat{\Lambda}_{i}\right\}  _{i=1}^{r}$ and $\left\{  \beta
_{i}^{n-3}\right\}  _{i=1}^{r}$ are well-defined. By making appropriate
choices of local orientation and base paths, we may arrange that
$\varepsilon_{\mathbb{Z[\pi}_{1}\left(  B_{2}\right)  /L]}\left(
\widehat{\Lambda}_{i},\beta_{j}^{n-3}\right)  =\pm\delta_{ij}$ for all $1\leq
i,j\leq r$. By applying part 2) of Lemma \ref{Lemma: Surgery on surfaces},
these surfaces may be surgered into a collection of $2$-spheres $\left\{
\Gamma_{i}\right\}  _{i=1}^{r}$ which is algebraically dual over
$\mathbb{Z[\pi}_{1}\left(  B_{2}\right)  /L]$ to the collection $\{\beta
_{i}^{n-3}\}$. Standard bundle theory again ensures that the $\Gamma_{i}$
inherit trivial normal bundles from the $\widehat{\Lambda}_{i}$. Attach
$3$-handles along regular neighborhoods of these $2$-spheres to obtain a
cobordism $\left(  W,A,B\right)  $. Then $A\hookrightarrow W$ induces a
$\pi_{1}$-isomorphism, with $\pi_{1}\left(  A\right)  \cong\pi_{1}\left(
W\right)  \cong\pi_{1}\left(  B_{2}\right)  $ and the above intersection data
assures that $A\hookrightarrow W$ is a simple $\mathbb{Z[\pi}_{1}\left(
A\right)  /L]$-equivalence.\ The surjectivity of $\pi_{1}\left(  B\right)
\rightarrow\pi_{1}\left(  W\right)  $ is clear from the construction.
\end{proof}

\begin{proof}
[Sketch of the Generalized Embedded Manifold Plus Construction]For the most
part, this proof follows the same outline as the proof of Theorem
\ref{Th: EPC} with modifications analogous to those found in the above sketch.
A few items become more delicate; we focus our attention on those issues.

1) In the proof of Theorem \ref{Th: EPC} we carried out the entire
construction inside a regular neighborhood $R^{\prime}$ of $T$, chosen early
in the proof. That was done solely for use in Step III. There it was crucial
that $B_{2}^{\left(  r\right)  }\hookrightarrow\overline{R^{\prime}-T^{\left(
r\right)  }}$ induce a $\pi_{1}$-surjection---thereby allowing us to choose an
arc $\alpha^{\prime\prime}$ in $B_{2}^{\left(  r\right)  }$ so that the loop
$\alpha\cup a^{\prime}\cup\alpha^{\prime\prime}$ contracted in $\overline
{R^{\prime}-T^{\left(  r\right)  }}$. In the more general case at hand, it
will be impossible to carry out the entire construction in a regular
neighborhood of $T$. Instead, we will expand the region where we work to an
open set $R^{\prime\prime}\supseteq T$ in which all loops corresponding to
elements of $L^{\prime}$ contract and for which $B_{2}\hookrightarrow
\overline{R^{\prime\prime}-\operatorname*{int}T}$ induces a $\pi_{1}%
$-surjection. If $L^{\prime}$ is the normal closure in $\pi_{1}\left(
B,\ast\right)  $ of a finite set of its elements, this is easy---let $S\approx
B\times\left[  0,1\right]  $ be a collar neighborhood of $B$ in $R$ with
$B_{1}$ the interior boundary component. Choose a pairwise disjoint collection
of properly embedded $2$-disks $\left\{  D_{1},\cdots,D_{r}\right\}  $ in
$\overline{R-S}$ whose boundaries in $B_{1}$ represent a finite normal
generating set for $G$, then supplement that collection with a disjoint
collection of pairwise disjoint $2$-disks $\left\{  D_{r+1},\cdots
,D_{s}\right\}  $ in $\overline{R-S}$ whose boundaries, together with those of
$\left\{  D_{1},\cdots,D_{r}\right\}  $, form a normal generating set for
$L^{\prime}$. Then $T$ may be viewed as a regular neighborhood of $S\cup
(\cup_{i=1}^{r}D_{i})$ in $R$ and we may let $R^{\prime\prime}$ be a regular
neighborhood of $S\cup(\cup_{i=1}^{s}D_{i})$ chosen to contain $T$ in its interior.

When $L^{\prime}$ is not normally finitely generated we use a similar, but
more delicate construction. Choose an infinite collection of $2$-disks
$\left\{  D_{r+1},D_{r+2},\cdots\right\}  $ whose boundaries, together with
those of $\left\{  D_{1},\cdots,D_{r}\right\}  $ generate $L^{\prime}$. These
may be chosen inductively so that each $D_{i}$ has a neighborhood $U_{i}$ for
which the collection $\{U_{i}\}_{i=1}^{\infty}$ is pairwise disjoint. We may
then thicken each $D_{i}$ ($i>r$) to a $2$-handle $h_{i}^{2}\subseteq U_{i}$
and add to a slightly enlarged copy of $T$ the interiors of each of these
$2$-handles. This creates an open subset $R^{\prime\prime}$ of $R$ containing
$T$ and having the desired properties.

2) Following the same strategy sketched out in the proof of Theorem
\ref{Th: GMPC}, but utilizing the more delicate item 2) of Lemma
\ref{Lemma: Surgery on surfaces}, we obtain a collection of $2$-spheres
$\left\{  \Gamma_{i}\right\}  _{i=1}^{r}$ in $B_{2}$ which is algebraically
dual over $\mathbb{Z[\pi}_{1}\left(  B_{2}\right)  /L]$ to the collection
$\{\beta_{i}^{n-3}\}_{i=1}^{r}$. (Note that $\omega_{1}\left(  L^{\prime
}\right)  \equiv1$ since $L^{\prime}\leq\ker\left(  \pi_{1}\left(
B,\ast\right)  \overset{i_{\ast}}{\longrightarrow}\pi_{1}\left(
R,\ast\right)  \right)  $.) Next we proceed inductively through the spherical
alteration process in the same manner as Step II of the proof of Theorem
\ref{Th: EPC} so that, at the conclusion, we have a new set of $2$-handles
$\left\{  \dot{h}_{i}\right\}  _{i=1}^{r}$ in $R^{\prime\prime}$ attached to
$S$ and a collection $\left\{  \dot{\Gamma}_{i}^{2}\right\}  _{i=1}^{r}$of
$2$-spheres algebraically dual over $\mathbb{Z[\pi}_{1}\left(  B_{2}^{\left(
r\right)  },\ast\right)  /L]$ to the collection $\left\{  \beta_{i}%
^{n-3}\right\}  _{i=1}^{r}$ of belt spheres of those $2$-handles in
$B_{2}^{\left(  r\right)  }$. In addition, each $\dot{\Gamma}_{i}^{2}$
contracts in $R^{\prime\prime}$. Contractibility of the $2$-spheres is more
delicate in this generalized situation. We use the full strength of Assertion
2) of Lemma \ref{alteration lemma}, the key point being that $j_{\#}%
:\mathbb{\pi}_{1}\left(  B_{2},\ast\right)  \rightarrow\mathbb{\pi}_{1}\left(
R^{\prime\prime},\ast\right)  $ is precisely the homomorphism that kills
$L\trianglelefteq\mathbb{\pi}_{1}\left(  B_{2},\ast\right)  $. By general
position these $2$-spheres also contract in $\overline{R^{^{\prime\prime}%
}-T^{\left(  r\right)  }}$, so for $n\geq7$ we may thicken a corresponding
collection of pairwise disjoint $3$-disks to $3$-handles to complete the
construction of $W$.

3) For $n=6$, Step III of the proof Theorem \ref{Th: EPC} goes through without
any changes. It is here, however, where we use the carefully chosen set
$R^{\prime\prime}$ in which to carry out the construction.
\end{proof}

\section{A more general lemma}

The following technical lemma was specifically designed for use in \cite{GT3}.
It is more general than Theorem \ref{Th: GEPC}, but no new ideas or techniques
are needed. For the reader who has made it this far, the proof is almost immediate.

\begin{lemma}
Let $R^{\prime}\subseteq R$ be a pair of $n$-manifolds ($n\geq6$) with a
common boundary component $B$, and suppose there is a subgroup $L^{\prime}$ of
$\ker\left(  \pi_{1}\left(  B\right)  \rightarrow\pi_{1}\left(  R\right)
\right)  $ for which $K=\ker\left(  \pi_{1}\left(  B\right)  \rightarrow
\pi_{1}\left(  R^{\prime}\right)  \right)  $ is strongly $L^{\prime}$-perfect.
Suppose further that there is a clean submanifold $T\subseteq R^{\prime}$
consisting of a finite collection $\mathcal{H}^{2}$ of $2$-handles in
$R^{\prime}$ attached to a collar neighborhood $S$ of $B$ with
$T\hookrightarrow R^{\prime}$ inducing a $\pi_{1}$-isomorphism (the
$2$-handles precisely kill the group $K$) and a finite collection $\left\{
\Theta_{t}^{2}\right\}  $ of pairwise disjoint embedded $2$-spheres in
$\partial T-B$, each of which contracts in $R^{\prime}$. Then on any
subcollection $\left\{  h_{j}^{2}\right\}  _{j=1}^{k}\subseteq\mathcal{H}^{2}%
$, one may perform spherical alterations to obtain to obtain $2$-handles
$\left\{  \dot{h}_{j}^{2}\right\}  _{j=1}^{k}$ in $R^{\prime}$ so that in
$\partial\dot{T}-B$ (where $\dot{T}$ is the correspondingly altered version of
$T$) there is a collection of $2$-spheres $\left\{  \dot{\Gamma}_{j}%
^{2}\right\}  _{j=1}^{k}$ algebraically dual over $\mathbb{Z[}\pi_{1}\left(
B\right)  /L^{\prime}]$ to the belt spheres $\left\{  \beta_{j}^{n-3}\right\}
_{j=1}^{k}$ common to $\left\{  h_{j}^{2}\right\}  _{j=1}^{k}$ and $\left\{
\dot{h}_{j}^{2}\right\}  _{j=1}^{k}$ with the property that each $\dot{\Gamma
}_{j}^{2}$ contracts in $R$. Furthermore, each correspondingly altered
$2$-sphere $\dot{\Theta}_{t}^{2}$ (now lying in $\partial\dot{T}-B$) has the
same $\mathbb{Z[}\pi_{1}\left(  B\right)  /L^{\prime}]$-intersection number
with those belt spheres and with any other oriented $\left(  n-3\right)
$-manifold lying in both $\partial T-B$ and $\partial\dot{T}-B$ as did
$\Theta_{t}^{2}$. Whereas the $2$-spheres $\left\{  \Theta_{t}^{2}\right\}  $
each contracted in $R^{\prime}$, the $\left\{  \dot{\Theta}_{t}^{2}\right\}  $
each contract in $R$.
\end{lemma}

\begin{remark}
\emph{In the usual way, when }$n\geq7$\emph{, the above result together with
general position assures the existence of a pairwise disjoint collection of
embedded }$3$\emph{-disks in }$\overline{R-\dot{T}}$\emph{ with boundaries
corresponding to the }$2$\emph{-spheres }$\left\{  \dot{\Gamma}_{j}%
^{2}\right\}  _{j=1}^{k}\cup\left\{  \dot{\Theta}_{t}^{2}\right\}  $\emph{.
Those }$3$\emph{-disks may be thickened to }$3$\emph{-handles by taking
regular neighborhoods. When }$n=6$\emph{, the same is true, but the }$\pi
$\emph{-}$\pi$\emph{ argument used in the proofs of Theorems \ref{Th: EPC} and
\ref{Th: GEPC} must again be employed to obtain embedded and pairwise disjoint
}$3$\emph{-disks. In that case we should also use the strategy used in the
proof of Theorem \ref{Th: GEPC} to ensure that we are working within a
submanifold }$R^{\prime\prime}$\emph{ of }$R$\emph{ which contains }$T$\emph{
and in which the group }$L^{\prime}$\emph{ dies and }$\pi_{1}\left(  B\right)
\rightarrow\pi_{1}\left(  R^{\prime\prime}\right)  $\emph{ is surjective.}
\end{remark}

\end{document}